\documentclass[11pt]{article}

\usepackage[utf8]{inputenc}
\usepackage{amsmath,amssymb,epsfig,bbm}
\usepackage{comment}
\usepackage[dvips]{color}

\usepackage{todonotes}
\usepackage{enumitem}
\setlist{nolistsep}

\newtheorem{defi}{Definition}
\newtheorem{prop}[defi]{Proposition}
\newtheorem{theo}[defi]{Theorem}
\newtheorem{conj}[defi]{Conjecture}
\newtheorem{lemm}[defi]{Lemma}
\newtheorem{coro}[defi]{Corollary}
\newtheorem{rema}[defi]{Remark}
\newtheorem{exem}[defi]{Example}
\newtheorem{exems}[defi]{Examples}

\newcommand{\bdefi}{\begin{defi}}
\newcommand{\edefi}{\end{defi}}
\newcommand{\bprop}{\begin{prop}}
\newcommand{\eprop}{\end{prop}}
\newcommand{\btheo}{\begin{theo}}
\newcommand{\etheo}{\end{theo}}
\newcommand{\blemm}{\begin{lemm}}
\newcommand{\brema}{\begin{rema}}
\newcommand{\erema}{\end{rema}}
\newcommand{\bexer}{\begin{exem}}
\newcommand{\eexer}{\end{exem}}
\newcommand{\bexems}{\begin{exems}}
\newcommand{\eexems}{\end{exems}}
\newcommand{\bconj}{\begin{conj}}
\newcommand{\econj}{\end{conj}}
\newcommand{\elemm}{\end{lemm}}
\newcommand{\bcoro}{\begin{coro}}
\newcommand{\ecoro}{\end{coro}}
\newcommand{\dem}{\noindent{\bf Proof. }}

\usepackage{mathrsfs}
\renewcommand\mathcal{\mathscr}

\newcommand{\A}{{\cal A}}

\newcommand{\M}{{\cal M}}
\newcommand{\N}{{\cal N}}

\newcommand{\D}{{\cal D}}

\newcommand{\V}{{\cal V}}
\newcommand{\W}{{\cal W}}

\renewcommand{\H}{{\cal H}}

\newcommand{\C}{{\cal C}}

\newcommand{\U}{{\cal U}}

\newcommand{\maths}[1]{{\mathbb #1}}  

\newcommand{\OO}{\maths{O}}
\newcommand{\RR}{\maths{R}}
\newcommand{\NN}{\maths{N}}
\newcommand{\CC}{\maths{C}}

\newcommand{\SSS}{\maths{S}}

\newcommand{\HH}{\maths{H}}
\newcommand{\FF}{\maths{F}}
\newcommand{\KK}{\maths{K}}
\newcommand{\ZZ}{\maths{Z}}

\newcommand{\ra}{\rightarrow}
\newcommand{\bs}{\backslash}

\newcommand{\ov}[1]{{\overline #1}} 
\newcommand{\wt}[1]{{\widetilde{#1}}}
\newcommand{\wh}[1]{{\widehat{#1}}}

\newcommand{\ga}{\gamma}
\newcommand{\Ga}{\Gamma}

\newcommand{\cqfd}{\hfill$\Box$}

\newcommand{\card}{{\operatorname{Card}}}

\newcommand{\Vol}{\operatorname{Vol}}

\newcommand{\PSL}{\operatorname{PSL}}

\newcommand{\diam}{\operatorname{diam}}
\newcommand{\Heis}{\operatorname{Heis}}

\newcommand{\hnr}{{\HH}^n_\RR}

\newcommand{\Convhull}{\C} 
\newcommand\mBM{m_{\rm BM}}
\newcommand\wtmBM{\wt m_{\rm BM}}

\newcommand{\stab}{\operatorname{Stab}}
\newcommand{\flow}[1]{{g^{#1}}}  
\newcommand\mussu[1]{\mu_{W^\pm(#1)}}

\newcommand\muss[1]{\mu_{W^+(#1)}}
\newcommand\musu[1]{\mu_{W^-(#1)}}

\newcommand{\normal}[1]{\partial^1_{+}{#1}}

\newcommand\Perp{\operatorname{Perp}}

\newcommand\normalout{\partial^1_{+}}
\newcommand\normalin{\partial^1_{-}}
\newcommand\normalpm{\partial^1_{\pm}}
\newcommand\normalmp{\partial^1_{\mp}}

\usepackage{enumitem}
\setlist[itemize]{labelindent=*}
\setlist[enumerate]{labelindent=*,leftmargin=*}
\usepackage{hyperref}

\newcounter{fig}

\title{Counting common perpendicular arcs in negative curvature}
\author{Jouni Parkkonen \and Fr\'ed\'eric Paulin} 

\begin{document}
\bibliographystyle{../alphanum}

\maketitle

\begin{abstract} 
  Let $D^-$ and $D^+$ be properly immersed closed locally convex subsets of a
  Riemannian manifold with pinched negative sectional curvature. Using
  mixing properties of the geodesic flow, we give an asymptotic
  formula as $t\ra+\infty$ for the number of common perpendiculars of
  length at most $t$ from $D^-$ to $D^+$, counted with multiplicities,
  and we prove the equidistribution in the outer and inner unit normal
  bundles of $D^-$ and $D^+$ of the tangent vectors at the endpoints
  of the common perpendiculars.  When the manifold is compact with
  exponential decay of correlations or arithmetic with finite volume,
  we give an error term for the asymptotic. As an application, we give
  an asymptotic formula for the number of connected components of the
  domain of discontinuity of Kleinian groups as their diameter goes to
  $0$.  
\end{abstract}

\section{Introduction}

Let $M$ be a complete connected Riemannian manifold with pinched
sectional curvature at most $-1$ whose fundamental group is not
virtually nilpotent, let $(\flow t)_{t\in\RR}$ be its geodesic
flow. Let $D^-$ and $D^+$ be proper nonempty properly immersed closed
locally convex subsets of $M$. A {\it common perpendicular} from $D^-$
to $D^+$ is a locally geodesic path in $M$ starting perpendicularly
from $D^-$ and arriving perpendicularly to $D^+$ (see Section
\ref{subsec:creatcomperp} for explanations when the boundary of $D^-$
or $D^+$ is not smooth.  For all $t> 0$, we denote by $\Perp (D^-,
D^+, t)$ the set of common perpendiculars from $D^-$ to $D^+$ with
length at most $t$ (considered with multiplicities), and by
$\N_{D^-,\,D^+}(t)$ its cardinality.  We refer to Section
\ref{subsec:downstairs} for the definition of the multiplicities,
which are equal to $1$ if $D^-$ and $D^+$ are embedded and disjoint.

In this paper, we give a general asymptotic formula $\N_{D^-,\,D^+}(t)
\sim c\,e^{c't}$ as $t\ra+\infty$, with error term estimates, and we
prove the equidistribution of the initial and terminal tangent vectors
of the common perpendiculars in the outer and inner unit normal
bundles of $D^-$ and $D^+$, respectively. The constants $c,c'$ are
explicit in terms of the Bowen-Margulis measure $\mBM$ of $M$ and the
skinning measures $\sigma^\mp_{D^\pm}$ of $D^\pm$.  These measures are
appropriate pushforwards of the Patterson-Sullivan densities to the
unit normal bundles of the lifts of $D^-$ and $D^+$ in the universal
cover of $M$ described in the present generality for the outer normal
bundle in \cite{ParPau14ETDS}, generalising \cite{OhSha12,OhSha13}
where $M$ has constant curvature and $D^-,D^+$ are balls, horoballs or
totally geodesic submanifolds.

We now state our counting and equidistribution results. We avoid any
compactness assumption on $M$, we only assume that the Bowen-Margulis
measure of $M$ is finite and that it is mixing for the geodesic
flow. We denote the total mass of any measure $m$ by $\|m\|$. Let
$\delta$ be the critical exponent of the fundamental group of $M$.

\btheo\label{theo:mainintrocount} Assume that the skinning measures
$\sigma^+_{D^-}$ and $\sigma^{-}_{D^+}$ are finite and nonzero. Then,
as $s\to+\infty$,
$$
\N_{D^-,\,D^+}(s)\sim
\frac{\|\sigma^+_{D^-}\|\,\|\sigma^{-}_{D^+}\|}
{\|\mBM\|}\,\frac{e^{\delta \,s}}{\delta}\,.
$$
\etheo

We refer to \cite{DalOtaPei00} for a finiteness criteria of $\mBM$, to
\cite{ParPau14ETDS} for finiteness criterion of the skinning measures,
generalising \cite{OhSha13}, and to \cite{Babillot02b} for mixing
criteria of $\mBM$.

The counting function $\N_{D^-,\,D^+}$ has been studied for particular
 triples $(M,D^-,D^+)$ at least since the 1950's for example
in \cite{Huber59}, \cite{Herrmann62}, \cite{Margulis69},
\cite{EskMcMul93}, \cite{Cosentino99}, 
\cite{Roblin03}, \cite{HerPau04}, \cite{Kontorovich09},
\cite{KonOh11}, \cite{OhSha12}, \cite{ParPau12JMD}, \cite{Kim13},
\cite{Pollicott11} and \cite{OhShaCircles}.  See 
\cite{ParPauRev} for a more detailed review.

When $M$ is a finite volume hyperbolic manifold, we get very explicit
forms of the counting results also in cases that were not known
before, see Corollary \ref{coro:constcurvcount}.  For example, if
$D^-$ and $D^+$ are closed geodesics of $M$ of lengths $\ell_-$ and
$\ell_+$, respectively, then the number of common perpendiculars from
$D^-$ to $D^+$ of length at most $s$ satisfies, as $s\ra+\infty$,
\begin{equation}\label{eqcaseclosedgeod}
\N_{D^-,\,D^+}(s)\sim
\frac{\pi^{\frac{n}{2}-1}\Ga(\frac{n-1}{2})^2}
{2^{n-2}(n-1)\Ga(\frac{n}{2})}
\;\frac{\ell_-\ell_+}{\Vol(M)}\;
e^{(n-1)s}\;.
\end{equation}
When $M$ is a closed surface and $D^-=D^+$, the formula
\eqref{eqcaseclosedgeod} is proved in \cite{MarMcKWam11} by trace
formula methods, but obtaining the case $D^-\neq D^+$ seems difficult
by their methods.

We denote the initial and terminal unit tangent vectors of $\alpha\in
\Perp (D^-, D^+, t)$ by $v^-_\alpha$ and $v^+_\alpha$, and the unit
Dirac mass at a point $z$ by $\Delta_z$. Theorem
\ref{theo:mainintrocount} is deduced in Section
\ref{sec:perpendiculars} from the following joint equidistribution
result of these vectors towards the skinning measures of $D^-$ and
$D^+$, that generalises work of Herrmann and Roblin for special
$D^\pm$'s.

\btheo \label{theo:mainintroequidis} For the weak-star convergence of
measures on $T^1M\times T^1M$, we have
$$
\lim_{t\ra+\infty}\; \delta\;\|\mBM\|\;e^{-\delta\, t}
\sum_{\alpha\in\Perp(D^-,\,D^+,\,t)} \; 
\Delta_{v^-_\alpha} \otimes\Delta_{v^+_\alpha}\;=\;
\sigma^+_{D^-}\otimes \sigma^-_{D^+}\,.
$$
\etheo

Both results are valid when $M$ is a good orbifold instead of a
manifold (for the appropriate notion of multiplicities), and when
$D^-,D^+$ are replaced by locally finite families.

Besides giving a unified treatment that covers all the special cases
in the above references, we have very weak finiteness and curvature
assumptions on the manifold, and no regularity assumptions on the
convex sets (see Corollary \ref{coro:geneOhShaKimintro} for a striking
application using convex sets with fractal boundary).  We develop new
techniques (necessary for the generality considered in
this paper) in order to apply Margulis's idea to use the mixing of the
geodesic flow: Due to the symmetry of the problem, a one-sided pushing of
the geodesic flow is not sufficient and we push simultaneously the
outer/inner unit normal vectors to the convex sets in opposite
directions. We also need a new effective study of the geometry and the
dynamics of the creation of common perpendiculars, see Subsection
\ref{subsec:creatcomperp}.

In the cases when the geodesic flow is known to be exponentially
mixing on $T^1M$, we obtain an exponentially small error term in the
approximation of the counting function $\N_{D^-,\,D^+}$. In particular,
when $M$ is arithmetic, the error term in Equation
\eqref{eqcaseclosedgeod} is $\operatorname{O}(e^{(n-1-\kappa)t})$ for
some $\kappa>0$.

\btheo\label{theo:4} Assume that $M$ is compact and the geodes\-ic
flow is exponentially mixing for the Hölder regularity, or that $M$ is
locally symmetric, the boundary of $D^\pm$ is smooth, $\mBM$ is
finite, smooth, and exponentially mixing under the geodesic flow for
the Sobolev regularity. Assume that the strong stable/unstable ball
masses by the conditionals of $\mBM$ are Hölder-continuous in their
radius.  Then there exists $\kappa>0$ such that, as $t\ra+\infty$,
$$
\N_{D^-,\,D^+,\, F}(t)=
\frac{\|\sigma^+_{D^-}\|\;\|\sigma^-_{D^+}\|}{\delta\;\|\mBM\|}\;
e^{\delta\, t}\big(1+\operatorname{O}(e^{-\kappa t})\big)\;.
$$
\etheo 

See Section \ref{sect:erroterms} for a discussion of the assumptions
and the dependence of $\operatorname{O}(\cdot)$ on the data. Similar
(sometimes more precise) error estimates were known earlier for the
counting function in special cases of $D^\pm$ in constant curvature
geometrically finite manifolds (often in small dimension) through the
work of Huber, Selberg, Patterson, Lax-Phillips \cite{LaxPhi82},
Cosentino \cite{Cosentino99}, Kontorovich-Oh \cite{KonOh11}, Lee-Oh
\cite{LeeOh13}.

\begin{center}
\includegraphics[width=120mm]{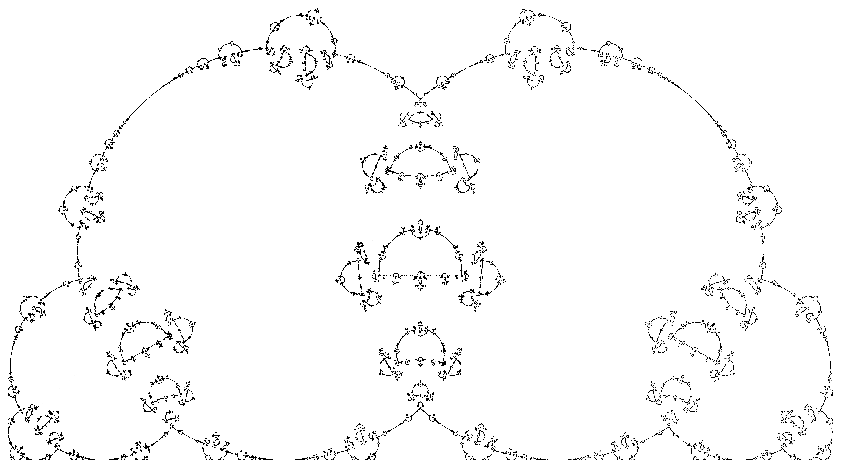}
\end{center}

Consider the picture above produced by D.~Wright's program
{\tt kleinian}, which is the limit set of a free product
$\Ga=\Ga_0*\ga_0\Ga_0\ga_0^{-1}$ of a quasifuchsian group $\Ga_0$ and
its conjugate by a large power $\ga_0$ of a loxodromic element whose
attracting fixed point is contained in the bounded component of
$\CC-\Lambda\Ga$, so that the limit set of $\Ga$ is the closure of a
countable union of quasi-circles.  As we will see in Section
\ref{subsec:countinlimset}, the number of Jordan curves in
$\Lambda\Ga$ with diameter at least $1/T$ is equivalent to
$c\,T^\delta$ where $c>0$ and $\delta\in \,]1,2[$ is the Hausdorff
dimension of the limit set.

\bcoro \label{coro:geneOhShaKimintro} Let $\Ga$ be a geometrically
finite discrete group of isometries of the upper halfspace model of
$\HH^n_\RR$, with bounded limit set $\Lambda\Ga$ in $\RR^{n-1}
= \partial_\infty \HH^n_\RR -\{\infty\}$.  Let $\delta$ be the
Hausdorff dimension of $\Lambda\Ga$.  Let $\Ga_0$ be a
convex-cocompact subgroup of $\Ga$ with infinite index.  Then, there
exists an explicitable $c>0$ such that, as $T\ra+\infty$,
$$
\card\{\ga\in \Ga/\Ga_0\;:\; \operatorname{diam}(\ga\Lambda\Ga_0)
\geq 1/T\}\sim \;c\,T^{\delta}\;.
$$
\ecoro

This corollary is due to \cite{OhShaCircles} when the limit set of
$\Ga_0$ is a round sphere (allowing the use of homogeneous
dynamics). We refer to Corollary \ref{coro:geneOhShaKimintrogene} for
a more general version and to Section \ref{subsec:countinlimset} for
complements and for extensions to any rank one symmetric space.

The results of this paper have been announced in the survey
\cite{ParPauRev}.  In \cite{ParPau14AFST}, we give several arithmetic
applications of the results of this paper, obtained by considering
arithmetically defined manifolds and orbifolds of constant negative
curvature.  In \cite{ParPauHeis}, we consider arithmetic applications
in the Heisenberg group via complex hyperbolic geometry.  In
\cite{ParPau15MZ}, we study counting and equidistribution in conjugacy
classes, giving a new proof of the main result of \cite{Huber56} and
generalising it to parabolic cyclic and more general subgroups,
arbitrary dimension, infinite volume and variable curvature.

The previous ArXiv version of this article contained the extension of
the counting and equidistribution results \ref{theo:mainintrocount},
\ref{theo:mainintroequidis} and \ref{theo:4} to Gibbs measures (that
is, equilibrium states associated with H\"older potentials on $T^1M$)
and counting functions with weights. In order to shorten this paper,
this material will appear as part of \cite{BroParPau16}.

\medskip\noindent{\small {\it Acknowledgement: } The first author
  thanks the Université Paris-Sud  for a month of visiting
  professor where this work was started, and the FIM of ETH Z\"urich
  for its support when this work was continued.  The
  second author thanks  ETH Z\"urich for frequent stays during
  the writing of this paper. Both authors thank the Mittag-Leffler
  Institute where this paper was almost completed. }

\section{Geometry, dynamics and convexity}
\label{sec:geomdyn}

Let $\wt M$ be a complete simply connected Riemannian manifold with
(dimension at least $2$ and) pinched negative sectional curvature
$-b^2\le K\le -1$, and let $x_0\in\wt M$ be a fixed basepoint.  Let
$\Ga$ be a nonelementary (not virtually nilpotent) discrete group of
isometries of $\wt M$, and let $M$ be the quotient Riemannian orbifold
$\Ga\bs\wt M$.  We denote by $\partial_{\infty}\wt M$ the boundary at
infinity of $\wt M$, by $\delta= \lim_{t\ra+\infty}\frac{1}{t}
\,\card\,\{\ga\in\Ga\,:\, d(x_0,\ga x_0)\leq t\}$ the critical
exponent of $\Ga$, and by $\Lambda\Ga$ the limit set of $\Ga$.

In this section, we review the required background on negatively
curved Riemannian manifolds seen as locally CAT($-\kappa$) spaces (see
\cite{BriHae99} for definitions, proofs and complements). We introduce
the notation for the outward and inward pointing unit normal bundles
of the boundary of a convex subset, and we define dynamical
thickenings in the unit tangent bundle of subsets of these
submanifolds, expanding on \cite{ParPau14ETDS}. We give a precise
definition of common perpendiculars in Subsection
\ref{subsec:creatcomperp} and we also give a procedure to construct
them by dynamical means.

For every $\epsilon>0$, we denote by $\N_\epsilon A$ the closed
$\epsilon$-neighbourhood of a subset $A$ of any metric space, by
$\N_{-\epsilon} A$ the set of points $x\in A$ at distance at least
$\epsilon$ from the complement of $A$, and by convention
$\N_0A=\overline{A}$.

\subsection{Strong stable and unstable foliations, and 
Hamenstädt's  distances}
\label{subsec:strongham}

We identify the unit tangent bundle $T^1N$ (endowed with Sasaki's
Riemannian metric) of a complete
Riemannian manifold $N$ with the set of its locally geodesic lines
$\ell:\RR\to N$, by the inverse of the map sending a (locally)
geodesic line $\ell$ to its (unit) tangent vector $\dot{\ell}(0)$ at
time $t=0$.  We denote by $\pi:T^1 N\to N$ the {\em basepoint
  projection}, given by $\pi(\ell)=\ell(0)$.

The {\em geodesic flow} on $T^1N$ is the smooth one-parameter group of
diffeomorphisms $(\flow{t})_{t\in\RR}$ of $T^1\wt M$, where $\flow{t}
\ell\,(s)=\ell(s+t)$, for all $\ell\in T^1N$ and $s,t\in\RR$.  The
action of the isometry group of $N$ on $T^1N$ by postcomposition (that
is, by $(\ga,\ell) \mapsto \ga\circ\ell$) commutes with the geodesic
flow.

When $\Ga$ acts without fixed points on $\wt M$, we have an
identification $\Ga\bs T^1\wt M=T^1M$.  More generally, we denote by
$T^1M$ the quotient Riemannian orbifold $\Ga\bs T^1\wt M$.  We use the
notation $(\flow{t})_{t\in\RR}$ also for the (quotient) geodesic flow
on $T^1 M$.

For every $v\in T^1\wt M$, let $v_-\in\partial_\infty\wt M$ and
$v_+\in\partial_\infty\wt M$, respectively, be the endpoints at
$-\infty$ and $+\infty$ of the geodesic line defined by $v$.  Let
$\partial_\infty^2\wt M$ be the subset of $\partial_\infty\wt
M\times\partial_\infty\wt M$ which consists of pairs of distinct
points at infinity of $\wt M$.  {\em Hopf's parametrisation} of
$T^1\wt M$ is the homeomorphism which identifies $T^1\wt M$ with
$\partial_\infty^2\wt M\times\RR$, by the map $v\mapsto(v_-,v_+,t)$,
where $t$ is the signed distance of the closest point to $x_0$ on the
geodesic line defined by $v$ to $\pi(v)$.  We have $\flow{s} (v_-,
v_+,t) = (v_-,v_+,t+s)$ for all $s\in\RR$, and for all $\ga\in \Ga$,
we have $\ga(v_-,v_+,t)=(\ga v_-,\ga v_+, t+t_{\ga,\,v_-,\,v_+})$
where $t_{\ga,\,v_-,\,\,v_+}\in\RR$ depends only on $\ga,v_-,v_+$.

Let $\iota: T^1\wt M\to T^1\wt M$ be the  {\em
  antipodal (flip) map} of $T^1\wt M$ defined by $\iota v=-v$ or,
using geodesic lines, by $\iota v:t\mapsto v(-t)$. In Hopf's
parametrisation, the antipodal map is the map $(v_-,v_+,t)\mapsto
(v_+,v_-,-t)$. We denote the quotient map of $\iota$ again by $\iota:
T^1 M\to T^1 M$, and call it the {\em antipodal map} of $T^1 M$. We
have $\iota\circ \flow{t}= \flow{-t} \circ \iota$ for all $t\in\RR$.

The {\em strong stable/unstable manifold} of $v\in T^1\wt M$ is 
$$
W^{\pm}(v)=\{v'\in T^1\wt M\;:\;d(v(t),v'(t))\to
0 \textrm{ as } t\to\pm\infty\}\,.
$$
The union for $t\in\RR$ of the images under $\flow t$ of the strong
stable manifold of $v\in T^{1}\wt M$ is the {\em stable manifold} $
W^{0+}(v)=\bigcup_{t\in\RR}\flow t W^{+}(v)$ of $v$, which consists of
the elements $v'\in T^1\wt M$ with $v'_+=v_+$. Similarly,
$W^{0-}(v)=\bigcup_{t\in\RR}\flow t W^{-}(v)$, which consists of the
elements $v'\in T^1\wt M$ with $v'_-=v_-$, is the {\em unstable
  manifold} $W^{0-}(v)$ of $v$. The maps from $\RR\times W^{\pm} (v)$
to $W^{0\pm} (v)$ defined by $(s,v')\mapsto \flow sv'$ are smooth
diffeomorphisms.

The strong stable manifolds, stable manifolds, strong unstable
manifolds and unstable manifolds are the (smooth) leaves of
topological foliations that are invariant under the geodesic flow and
the group of isometries of $\wt M$, denoted by $\W^{+}, \W^{0+},
\W^{-}$ and $\W^{0-}$ respectively.  These foliations are
H\"older-continuous when $\wt M$ has compact quotients or when $\wt M$
has pinched negative sectional curvature with bounded derivatives (see
for instance \cite{Brin95}, \cite[Thm.~7.3]{PauPolSha}) and even
smooth when $\wt M$ is symmetric.

For any point $\xi\in\partial_{\infty}\wt M$, let $\rho_{\xi}:
[0,+\infty[\; \ra \wt M$ be the geodesic ray with origin $x_{0}$ and
point at infinity $\xi$.  The {\em Busemann cocycle} of $\wt M$ is the
map $\beta: \wt M\times\wt M\times\partial_{\infty} \wt M\to\RR$
defined by $(x,y,\xi)\mapsto \beta_{\xi}(x,y)=
\lim_{t\to+\infty}d(\rho_{\xi}(t),x)-d(\rho_{\xi}(t),y)\;.$

\smallskip\noindent \begin{minipage}{9.9cm} ~~~ The projections in
  $\wt M$ of the strong unstable and strong stable manifolds of $v\in
  T^1\wt M$, denoted by $H_{-}(v)=\pi(W^{-}(v))$ and $H_{+}(v)=
  \pi(W^{+}(v))$, are the {\em unstable and stable horospheres}
  of $v$ {\it centered at} $v_{-}$ and $v_{+}$, respectively. The
  unstable horosphere of $v$ coincides with the zero set of the map
  $x\mapsto f_{-}(x)= \beta_{v_-}(x,\pi(v))$ and the stable horosphere
  of $v$ is the zero set of $x\mapsto f_{+}(x)=
  \beta_{v_+}(x,\pi(v))$.  The sublevel sets $H\!B_{-}(v)=
  f_{-}^{-1}(]-\infty,0])$ and $H\!B_{+}(v)= f_{+}^{-1}(]-\infty,0])$
  are the {\em horoballs} bounded by $H_{-}(v)$ and
  $H_{+}(v)$. Horoballs are (strictly) convex subsets of $\wt
  M$. \end{minipage}
\begin{minipage}{5cm}
\begin{center}
\input{fig_horobull.pstex_t}
\end{center}
\end{minipage}

\smallskip For every $v\in T^1\wt M$, let $d_{W^{-}(v)}$ and
$d_{W^{+}(v)}$ be {\it Hamenst\"adt's distances} on the strong
unstable and strong stable leaf of $v$, defined as follows (see for
instance \cite[Appendix]{HerPau97}): for all $w,z\in W^{\mp}(v)$, let
$$
d_{W^{\mp}(v)}(w,z) = 
\lim_{t\ra+\infty} e^{\frac{1}{2}d(w(\pm t),\;z(\pm t))-t}\;.
$$
Hamenst\"adt's distances are distances inducing the original topology
on $W^{\pm}(v)$. For all $w,z\in W^{\pm}
(v)$ and  $s\in\RR$, and for every isometry $\ga$ of
$\wt M$, we have
\begin{equation}\label{eq:expandcontracHamdist}
d_{W^{\pm}(\ga v)}(\ga w,\ga z)= d_{W^{\pm}(v)}(w,z)\quad {\rm
and }\quad
d_{W^{\pm} (\flow sv)}(\flow sw,\flow sz)=e^{\mp s}d_{W^{\pm}(v)}(w,z)\;.
\end{equation}
A proof of the following result, first given in a preliminary version
of this paper, may now be found in \cite[Lem.~2.4]{PauPolSha}.

\blemm \label{lem:minodistHam} For all $v\in T^1\wt M$, $v'\in W^{\pm}
(v)$, we have
$$
d(\pi(v),\pi(v'))\leq d_{W^{\pm}(v)}(v,v')\;. \Box
$$
\elemm

\subsection{Dynamical thickening of outer and inner unit normal
  bundles}
\label{subsec:dynathick}

Let $D$ be a nonempty proper closed convex subset in $\wt M$. We
denote by $\partial D$ its boundary in $\wt M$ and by $\partial_\infty
D$ its set of points at infinity. In this subsection, we recall from
\cite{ParPau14ETDS} the definition of the outer unit normal bundle of
$\partial D$, the dynamical thickenings of its subsets, and we extend
these definitions to the inner unit normal bundle of $\partial D$.

Let $P_D:\wt M\cup(\partial_\infty\wt M-\partial_\infty D)\to D$ be
the (continuous) {\em closest point map} defined on
$\xi\in\partial_\infty\wt M -\partial_\infty D$ by setting $P_D(\xi)$
to be the unique point in $D$ that minimises the function $y\mapsto
\beta_\xi(y,x_0)$ from $D$ to $\RR$.  The {\em outer unit normal
  bundle} $\normalout D$ of the boundary of $D$ is the topological
submanifold of $T^1\wt M$ consisting of the geodesic lines $v:\RR\to
\wt M$ with $P_{D}(v_+) =v(0)$. The {\em inner unit normal bundle} of
the boundary of $D$ is $\normalin D=\iota\normalout D$. Note that
$\pi(\normalpm D)=\partial D$, that $\normalout H\!B_-(v)$ is the
strong unstable manifold $W^{-}(v)$ of $v$ and that $W^{+}(v)=
\normalin H\!B_+(v)$. When $D$ is a totally geodesic submanifold of
$\wt M$, then $\normalout D=\normalin D$.

The restriction of  $P_D$ to $\partial_\infty\wt
M-\partial_\infty D$ is not injective in general, but the inverse
$P^+_D$ of the restriction to $\partial^1_+D$ of the {\em (positive)
  endpoint map} $v\mapsto v_+$ is a natural lift of $P_D$ to a
homeomorphism from $\partial_\infty\wt
M- \partial_\infty D$ to $\normalout D$ such that $\pi\circ P^+_D =
P_D$.  Similarly, $P^-_{D}=\iota\circ P^+_D:\partial_\infty\wt
M-\partial_\infty D\to \normalin D$ is a 
homeomorphism such that $\pi\circ P^-_D=P_D$.

For every isometry $\ga$ of $\wt M$, we have $\normalpm (\ga D)=
\ga\,\normalpm D$ and $P^\pm_{\ga D}\circ\ga=\ga\circ P^\pm_D$. In
particular, $\normalpm D$ is invariant under the isometries of $\wt M$
that preserve $D$. For all $t\geq 0$, we have $\flow{\pm t}\normalpm
D= \normalpm (\N_tD)$.

We define 
\begin{equation*}\label{eq:defiUCp}
\U^\pm_D=\{v\in T^1\wt M:\ v_\pm\notin\partial_\infty D\}\;.
\end{equation*}
Note that $\U^-_D=\iota \U^+_D$, and that $\U^\pm_D$ is an open subset
of $T^1\wt M$, invariant under the geodesic flow.  We have
$\U^\pm_{\ga D}=\ga \U^\pm_D$ for every isometry $\ga$ of $\wt M$ and,
in particular, $\U^\pm_D$ is invariant under the isometries of $\wt M$
preserving $D$.

\medskip \noindent
\begin{minipage}{10.4cm} ~~~ Define a fibration $f^+_D:\U^+_D\ra
  \normalout D$ as the composition of the positive endpoint map from
  $\U^+_D$ onto $ \partial_\infty \wt M-\partial_\infty D$ (which is a
  fibration) and the homeomorphism $P^+_{D}$ from $\partial_\infty \wt
  M-\partial_\infty D$ to $\normalout D$.  The fiber of $w\in
  \normalout D$ for $f^+_D$ is exactly the stable leaf
  $$W^{0+}(w)=\{v\in T^1\wt M\;:\;v_+=w_+\}\,.$$
\end{minipage}
\begin{minipage}{4.5cm}
\begin{center}
\input{fig_UCfCnuPC.pstex_t}
\end{center}
\end{minipage}

\medskip
\noindent Analogously, we define a fibration $f^-_D= \iota\circ
f^+_D\circ \iota:\U^-_D\ra \normalin D$ as the composition of the
negative endpoint map with $P^-_D$, for which the fiber of
$w\in \normalin D$ is the unstable leaf $W^{0-}(w)=\{v\in T^1\wt
M\;:\;v_-=w_-\}.$

\smallskip For every isometry $\ga$ of $\wt M$, we have 
\begin{equation}\label{eq:equivfibrationf}
f^\pm_{\ga D} \circ\ga=\ga\circ f^\pm_D\,.
\end{equation}
We have $f^\pm_{\N_tD}= \flow{\pm t}\circ f^\pm_D$ for all $t\geq 0$,
and $f^\pm_D\circ \flow t=f^\pm_D$ for all $t\in\RR$. In particular,
the fibrations $f^\pm_D$ are invariant under the geodesic flow.

The next result will only be used for the error term estimates in
Section \ref{sect:erroterms}. Note that if $\wt M$ is a symmetric
space (in which case the strong stable and unstable foliations are
smooth, and the sphere at infinity has a smooth structure such that
the maps $v\mapsto v_\pm$ from $W^{\mp}(w)$ to $\partial_\infty\wt
M-\{w_\mp\}$ are smooth), and if $D$ has smooth boundary, then the
fibrations $f^\pm_D$ are smooth.

Recall that a map $f:X\ra Y$ between two metric spaces is (uniformly
locally) {\it H\"older-continuous} if there exist $c,c'>0$ and
$\alpha\in\;]0,1]$ such that $d(f(x),f(y))\leq c\;d(x,y)^\alpha$ for
all $x,y\in X$ with $d(x,y)\leq c'$.

\blemm \label{lem:fibrationholder}
The maps $f^\pm_D$ are Hölder-continuous on the set of elements $v\in
\U^\pm_D$ such that $d(\pi(v),\pi(f^\pm_D(v)))$ is bounded.  
\elemm

\dem We prove the result for $f^+_D$, the one for $f^-_D$ follows
similarly.  For all $u,u'\in T^1\wt M$, denote the geodesic lines they
define by $t\mapsto u_t,u'_t$, and let
$$
\delta_1(u,u')=\exp({-\sup\{t\geq 0\;:\;\sup_{s\in[-t,t]}d(u_s,u'_s)
  \leq 1\}})\;\;\;{\rm and}\;\;\; 
\delta_2(u,u')=\sup_{t\in[0,1]}d(u_t,u'_t)\;.
$$
with the convention $\delta_1(u,u')=1$ if $d(u_0,u'_0)> 1$ and
$\delta_1(u,u')=0$ if $u=u'$.  By for instance
\cite[p.~70]{Ballmann95}, the maps $\delta_1,\delta_2$ are distances
on $T^1\wt M$ which are Hölder-equivalent to Sasaki's distance.

\begin{center}
\input{fig_fibraholder.pstex_t}
\end{center}

Let $v,v'\in T^1\wt M$ be such that $d(v_0,v'_0)\leq 1$, let $w=f^+_D(v)$
and $w'=f^+_D(v')$. Let $T=\sup\{t\geq 0\;:\; \sup_{s\in[0,t]}
d(v_s,v'_s) \leq 1\}$, so that $\delta_1(v,v')\geq e^{-T}$. We may
assume that $T$ is finite, otherwise $v_+=v'_+$, hence $w=w'$. Let
$x=v_T$ and $x'=v'_T$, which satisfy $d(x,x')\leq 1$. Let $y$
(respectively $y'$) be the closest point to $x$ (respectively $x'$) on
the geodesic ray defined by $w$ (respectively $w'$). By convexity,
since $d(v_0,w_0)$ and  $d(v'_0,w'_0)$ are bounded by a constant $c>0$ and
since $v_+=w_+,v'_+=w'_+$, we have $d(x,y)\leq c$ and $d(x',y')\leq
c$. By the triangle inequality, we have $d(y,y')\leq 2c+1$,
$d(y,w_1)\geq T-2c-1$ and $d(y',w'_1)\geq T-2c-1$. By convexity, and
since projection maps exponentially decrease the distances, there
exists a constant $c'>0$ such that
$$
\delta_2(w,w')=d(w_1,w'_1)\leq c'd(y,y') e^{-(T-2c-1)}
\leq c'(2c+1)e^{2c+1}\;\delta_1(v,v')\,.
$$
The result follows.
\cqfd

\medskip 
\noindent
\begin{minipage}{10.4cm} ~~~ 
Let $\eta,\eta'>0$. For all $w\in T^1\wt M$, let
\begin{equation}\label{eq:defiboulehamen}
B^\pm(w,\eta')=\{v'\in W^{\pm}(w)\;:\;d_{W^{\pm}(w)}(v',w)<\eta'\}
\end{equation}
be the open balls of radius $\eta'$ centered at $w$ for Hamenst\"adt's
distance in the strong stable/unstable leaves of $w$.   Let
$$
V^\pm_{w,\,\eta,\,\eta'} =\bigcup_{s\in\;]-\eta,\,\eta\,[}\flow s
B^\pm(w, \eta')\;.
$$
\end{minipage}
\begin{minipage}{4.5cm}
\begin{center}
\input{fig_voisdyna.pstex_t}
\end{center}
\end{minipage}

\medskip
\noindent We have $B^-(w,\eta')=\iota B^+(\iota w,\eta')$ hence
$V^-_{w,\,\eta,\,\eta'}=\iota V^+_{\iota w,\,\eta,\,\eta'}$. We have 
\begin{equation}\label{eq:flotVetaeta}
\flow s B^\pm(w, \eta') = B^\pm(\flow sw, e^{\mp s} \eta')\;\; 
{\rm hence}\;\; \flow s
V^\pm_{w, \,\eta,\,\eta'} = V^\pm_{\flow sw, \,\eta,\,e^{\mp s} \eta'}
\end{equation} 
for all $s\in\RR$. For every isometry $\ga$ of $\wt M$, we have $\ga
B^\pm(w, \eta')= B^\pm(\ga w, \eta') $ and $\ga
V^\pm_{w,\,\eta,\,\eta'} = V^\pm_{\ga w,\,\eta,\,\eta'}\,$. The map
from $]-\eta,\eta[\;\times B^\pm(w, \eta')$ to
$V^\pm_{w,\,\eta,\,\eta'}$ defined by $(s,v') \mapsto \flow s v'$ is a
homeomorphism. For all subsets $\Omega^-$ of $\normalout D$ and
$\Omega^+$ of $\normalin D$, let
$$
\V^+_{\eta,\,\eta'}(\Omega^-)=
\bigcup_{w\in\Omega^-} V^+_{w,\,\eta,\,\eta'}
\;\;\;{\rm and}\;\;\;
\V^-_{\eta,\,\eta'}(\Omega^+)=
\bigcup_{w\in\Omega^+} V^-_{w,\,\eta,\,\eta'}\;.
$$
For every isometry $\ga$ of $\wt M$, we have $\ga\V^\pm_{\eta,\,\eta'}
(\Omega^{\mp})=\V^\pm_{\eta,\,\eta'}(\ga \Omega^{\mp})$ and for every
$t\geq 0$, we have
\begin{equation}\label{eq:flowbehavdynaneigh}
\flow {\pm t}\V^\pm_{\eta,\,\eta'}(\Omega^{\mp})
=\V^\pm_{\eta,\,e^{-t}\eta'}(\flow{\pm t}\Omega^{\mp})\;. 
\end{equation}
The thickenings (or {\it dynamical neighbourhoods})
$\V^\pm_{\eta,\,\eta'}(\Omega^{\mp})$ are nondecreasing in $\eta$ and
in $\eta'$ and their intersections and unions satisfy
$$
\bigcap_{\eta\,,\;\eta'>0}\V^\pm_{\eta,\,\eta'}(\Omega^{\mp})=\Omega^{\mp}
\;\;\;{\rm and}\;\;\;
\bigcup_{\eta\,,\;\eta'>0}\V^\pm_{\eta,\,\eta'}(\partial^1_{\pm} D)=\U^\pm_D\;.
$$ 
The restriction of $f^\pm_D$ to $\V^\pm_{\eta,\,\eta'}(\Omega^{\mp})$
is a fibration over $\Omega^{\mp}$, whose fiber over $w\in
\Omega^{\mp}$ is the open subset $V^\pm_{w,\,\eta,\,\eta'}$ of the
stable/unstable leaf of $w$.

\subsection{Creating common perpendiculars}
\label{subsec:creatcomperp}

For any two closed convex subsets $D^-$ and $D^+$ of $\wt M$, we say
that a geodesic arc $\alpha:[0,T]\to \wt M$, where $T>0$, is a {\em
  common perpendicular} from $D^-$ to $D^+$ if its initial tangent
vector $\dot\alpha(0)$ belongs to $\normalout D^-$ and if its terminal
tangent vector $\dot\alpha(T)$ belongs to $\normalin D^+$. It is
important to think of common perpendiculars as oriented arcs (from
$D^-$ to $D^+$). Note that there exists a common perpendicular from
$D^-$ to $D^+$ if and only if $D^-$ and $D^+$ are nonempty and the
closures $\overline{D^-}$ and $\overline{D^+}$ of ${D^-}$ and ${D^+}$
in the compactification $\wt M\cup\partial_\infty\wt M$ are
disjoint. A common perpendicular from $D^-$ to $D^+$,
if it exists, is unique.

When $\overline{D^-}$ and $\overline{D^+}$ are disjoint, and when
$\partial D^-$ and $\partial D^+$ are $\C^{1} $-submanifolds (for
instance, if $D^\pm$ are closed $\epsilon $-neighbourhoods of nonempty
convex subsets of $\wt M$ for some $\epsilon>0$, see \cite{Walter76}),
this definition of a common perpendicular corresponds to the usual
one. But there are interesting closed convex subsets with less regular
boundary, such as in general the convex hulls of limit sets of
nonelementary discrete groups of isometries of $\wt M$. Although it
would be possible to take the closed $\epsilon$-neighbourhood, to
count common perpendiculars in the usual sense, and then to take a
limit as $\epsilon$ goes to $0$, it is more natural to work directly
in the above generality (see \cite[Sect.~3.2]{ParPauRev} for further
comments).

The crucial observation is that two nonempty proper closed convex
subsets $D^-$ and $D^+$ of $\wt M$ have a common perpendicular
$\alpha$ of length a given $t> 0$ if and only if the pushforwards and
pullbacks by the geodesic flow at time $\frac{t}{2}$ of the outer and
inner normal bundles of $D^-$ and $D^+$, that is the subsets
$\flow{\frac{t}{2}}\normalout D^-$ and $\flow{-\frac{t}{2}} \normalin
D^+$ of $T^1\wt M$, intersect. Then their intersection is the
singleton consisting of the tangent vector of $\alpha$ at its
midpoint.

\blemm \label{lem:creationperp} For every $R>0$, there exist
$t_0,c_0>0$ such that for all $\eta\in\;]0,1]$ and all $t\in \;
[t_0,+\infty[\,$, for all nonempty closed convex subsets $D^-,D^+$ in
$\wt M$, and for all $w\in\flow{t/2}\V^+_{\eta,\,R}(\normalout D^-)
\cap \flow{-t/2} \V^-_{\eta,\,R}(\normalin D^+)$, there exist $s\in\;
]-2\eta,2\eta[$ and a common perpendicular $\wt c$ from $D^-$ to $D^+$
such that
\begin{itemize}
\item the length of $\wt c$ is contained in
  $[t+s-c_0\,e^{-\frac{t}{2}},t+s+c_0\,e^{-\frac{t}{2}}]$,
\item if $w^\mp=f^\pm_{D^\mp}(w)$ and if $p^\pm$ is the endpoint of
  $\wt c$ in $D^\pm$, then $d(\pi(w^\pm),p^\pm)\leq
  c_0\,e^{-\frac{t}{2}}$,
\item the basepoint $\pi(w)$ of $w$ is at distance at most
$c_0\,e^{-\frac{t}{2}}$ from a point of $\wt c$, and 
$$
\max\{\;d(\pi(g^{\frac{t}{2}} w^-),\pi(w)),\;
d(\pi(g^{-\frac{t}{2}} w^+),\pi(w))\;\}\leq 
\eta + c_0\,e^{-\frac{t}{2}}\;.
$$
\end{itemize}  
\elemm

\begin{center}
\input{fig_creationperp.pstex_t}
\end{center}

\dem Let $t\geq 3$ and $\eta\in\;]0, 1]$. By definition of the
dynamical neighbourhoods $\V^\mp_{\eta,\,R}(\normalmp D^\pm)$, there exist
$w^\pm\in\normalmp D^\pm$ and $s^\pm\in\;]-\eta,+\eta[$ such that
$$
d_{W^{+}(w^-)}(\flow{-\frac{t}{2}-s^-}w,w^-)\leq R\;\;\;{\rm  and}
\;\;\;d_{W^{-}(w^+)}(\flow{\frac{t}{2}+s^+}w,w^+)\leq R\;.
$$
Let $x^\pm=\pi(w^\pm)$, $y=\pi(w)$, and let $\alpha^-$ (respectively
$\alpha^+$) be the angle at $x^-$ (respectively $x^+$) between $w^-$
(respectively $\iota w^+$) and the geodesic segment $[x^-,x^+]$.

\medskip\noindent
\begin{minipage}{9cm}{\bf Step 1. } Let $\ov\alpha^-$ (respectively
  $\ov\alpha^+$) be the angle at $x^-$ (respectively $x^+$) between
  the outer normal vector $w^-$ (respectively $\iota w^+$) and the
  geodesic segment $[x^-,y]$ (respectively $[x^+,y]$).  Let
  $\beta^\pm$ be the angle at $y$ between $\pm w$ and the geodesic
  segment $[y,x^\pm]$. Let us prove that there exist two constants
  $t_1,c_1>0$ depending only on $R$ such that if $t\geq t_1$ then
  $\ov\alpha\,^\pm,\beta^\pm\leq c_1\,e^{-\frac{t}{2}}$.
\end{minipage} 
\begin{minipage}{5.9cm}
\begin{center}
\input{fig_angle1.pstex_t}
\end{center}
\end{minipage}

\medskip By Lemma \ref{lem:minodistHam} and Equation
\eqref{eq:expandcontracHamdist}, we have
$$
d(\pi(\flow{-\frac{t}{2}-s^-}w),x^-)\leq 
d_{W^{+}(w^-)}(\flow{-\frac{t}{2}-s^-}w,w^-)\leq 
R\;,
$$
\begin{equation}\label{eq:pincexpofortstab}
d(y,\pi(\flow{\frac{t}{2}+s^-}w^-))\leq
d_{W^{+}(w)}(w,\flow{\frac{t}{2}+s^-}w^-)
\leq R\,e^{-\frac{t}{2}-s^-}\;.
\end{equation}
In particular,
$$d(\pi(w),\pi(g^{\frac{t}{2}} w^-))\leq 
d(y,\pi(g^{\frac{t}{2}+s^-} w^-))+ |s^-|\leq 
R\,e^{-\frac{t}{2}-s^-}+\eta \leq c_0\;e^{-\frac{t}{2}}+\eta
$$
if we assume, as we may, that $c_0\geq R\,e^\eta$. With a similar
argument for $w^+$, this proves the last formula of Lemma
\ref{lem:creationperp}.

Recall that by a hyperbolic trigonometric formula (see for instance
\cite[p.~147]{Beardon83}), for any geodesic triangle in the real
hyperbolic plane, with angles $\alpha,\beta,\ga$ and opposite side
lengths $a,b,c$, if $\ga\geq \frac{\pi}{2}$, then $\tan \alpha\leq
\frac{1}{\sinh b}$, which is at most $\frac{1}{\sinh (c-a)}$ if $c> a$
by the triangle inequality. By comparison, if $t\geq 2(R+2)$ (which
implies that $\frac{t}{2}+s^--R\ge 1$), we hence have
$$
\max\{\tan \ov\alpha\,^-,\tan \beta^-\}\leq
\frac{1}{\sinh(\frac{t}{2}+s^--R)}\leq 4\,e^{-\frac{t}{2}-s^-+R}
$$ 
With a symmetric argument for $\ov\alpha\,^+,\beta^+$, the result
follows.

\bigskip\noindent
\begin{minipage}{8.5cm}{\bf Step 2. } Let $\overline{\ov\alpha}\,^\pm$
  be the angles at $x^\pm$ of the geodesic triangle with vertices
  $x_-,x_+,y$.  Let $y'$ be the closest point to $y$ on the side
  $[x_-,x_+]$. Let us prove that there exist two constants $t_2,c_2>0$
  depending only on $R$ such that if $t\geq t_2$ then
  $\overline{\ov\alpha}\,^\pm, d(y,y')\leq c_2\,e^{-\frac{t}{2}}$.
\end{minipage} 
\begin{minipage}{6.4cm}
\begin{center}
\input{fig_angle2.pstex_t}
\end{center}
\end{minipage} 

\medskip Since the angle $\angle_y(x^-,x^+)$ is at least $\pi-\beta^-
-\beta^+$, at least one of the two angles $\angle_y(y',x^\pm)$ is at
least $\frac{\pi-\beta^--\beta^+}{2}$. By a comparison argument
applied to one of the two triangles with vertices $(y,y',x^\pm)$, as
in the end of the first step, we have
$\tan \frac{\pi-\beta^--\beta^+}{2}\leq
\frac{1}{\sinh d(y,\,y')}$. Hence
$$
d(y,y')\leq \sinh \;d(y,y')\leq \tan \frac{\beta^-+\beta^+}{2}\;,
$$
and the desired majoration of $d(y,y')$ follows from Step 1.  By the
same argument, we have $\tan \overline{\ov\alpha}\,^\pm\leq
\frac{1}{\sinh (d(x^\pm,\,y)-d(y,\,y'))}$. Since $d(x^\pm,\,y)\geq
\frac{t}{2}-s^\pm-R\,e^{-\frac{t}{2}-s^-}$ by the inverse triangle
inequality and Equation \eqref{eq:pincexpofortstab}, the desired
majoration of $\overline{\ov\alpha}\,^\pm$ follows.

\medskip\noindent{\bf Step 3. } Let us prove that there exist two
constants $t_3,c_3>0$ depending only on $R$ such that if $t\geq t_3$
then there exists a common perpendicular $\wt c=[p^-,p^+]$ from $D^-$
to $D^+$ such that $d(x^-,p^-), d(x^+,p^+)\leq
c_3\,e^{-\frac{t}{2}}$. This will prove the second point of Lemma
\ref{lem:creationperp} (if $t_0\geq t_3$ and $c_0 \geq c_3$).

By the first two steps, we have, if $t\geq \min\{t_1,t_2\}$,
\begin{equation}\label{eq:alphapmmajo}
\alpha^\pm\leq \ov\alpha\,^\pm+\overline{\ov\alpha}\,^\pm\leq 
(c_1+c_2)e^{-\frac{t}{2}}\;.
\end{equation}
Assume by absurd that the intersection of the closures of $D^-$ and
$D^+$ in $\wt M\cup \partial_\infty\wt M$ contains a point $z$.  Then
by convexity of $D^\pm$, and since the distance $d(x^-,x^+)$ is large
and the angles $\alpha^\pm$ are small if $t$ is large, the angles at
$x^\pm$ of the geodesic triangle with vertices $z,x^-,x^+$ are almost
at least $\frac{\pi}{2}$, which is impossible since $\wt M$ is
$\operatorname{CAT}(-1)$.  Hence the nonempty closed convex subsets
$D^-$ and $D^+$ have a common perpendicular $\wt c=[p^-,p^+]$, with
$p^\pm\in D^\pm$.

Consider the geodesic quadrilateral $Q$ with vertices
$x^\pm,p^\pm$. By convexity of $D^\pm$, its angles at $p^\pm$ are at
least $\frac{\pi} {2}$ and its angles at $x^\pm$ are at least
$\frac{\pi} {2}- \alpha^\pm$.  Note that if $t\geq t'_2=2(R+c_2+1+
\operatorname{argsinh} 2)$ then we have, by Step 2,
\begin{align}
d(x^-,x^+)&\geq d(x^-,y)+d(y,x^+)-2\,d(y,y')\nonumber\\ & \geq 
(\frac{t}{2}+s^--R)+(\frac{t}{2}+s^+-R)-2\,c_2\,e^{-\frac{t}{2}}
\geq 2\operatorname{argsinh} 2\;.\label{eq:argsinhmino}
\end{align}

Up to replacing $Q$ by a comparison quadrilateral (obtained by gluing
two comparison triangles) in the real hyperbolic plane $\HH^2_\RR$,
having the same side lengths and bigger angles, we may assume that
$\wt M=\HH^2_\RR$ and that $x^-$ and $x^+$ are on the same side of the
geodesic line through $p^-,p^+$. Up to replacing $Q$ by a
quadrilateral having same distances $d(x^-,x^+)$, $d(x^-,p^-)$,
$d(x^+,p^+)$ and bigger angles at $x^-,x^+$, we may assume that the
angles at $p^-,p^+$ are exactly $\frac{\pi}{2}$. If, say, the angle at
$x^+$ was bigger than the angle at $x^-$, up to replacing $x^+$ by a
point on the geodesic line through $p^+,x^+$ on the other side of
$x^+$ than $p^+$ if $p^+\neq x^+$, which increases $d(x^-,x^+)$,
$d(x^+,p^+)$, decreases the angle at $x^+$ and increases the angle at
$x^-$, we may assume that the angles at $x^\pm$ are equal, and we
denote this common value by $\phi\geq
\frac{\pi}{2}-\min\{\alpha^-,\alpha^+\}$.

\medskip\noindent
\begin{minipage}{8.7cm}~~~ Let $b_1=\frac{1}{2}\,d(x^-,x^+)$ and
  $b_2=d(x^-,p^-)=d(x^+,p^+)$. By formulas of \cite[p.~157]{Beardon83}
  on Lambert quadrilaterals, we have
$$
\cosh b_2 =
\frac{\sinh b_{1}}{\sqrt{\sinh^2b_1\sin^2\phi-\cos^2\phi}}\;.
$$
\end{minipage} 
\begin{minipage}{6.2cm}
\begin{center}
\input{fig_quadranhyp.pstex_t}
\end{center}
\end{minipage}

\medskip \noindent By Equation \eqref{eq:alphapmmajo}, with
$c'_2=c_1+c_2$, let $t''_2>0$ be a constant, depending only on $R$,
such that if $t\geq t''_2$, then $\sin \phi\geq \max \{\cos
\alpha^\pm\}\geq 1-{c'_2}^2\,e^{-t}\geq 1/2$.  By Equation
\eqref{eq:argsinhmino}, if $t\geq t'_2$, then $b_1\geq \frac{t}{2}-R-1
-c_2\geq \operatorname{argsinh} 2$ (and in particular $1/\sinh b_1\leq
1/2$). Hence, if $t\geq\max\{t'_2,t''_2)\}$, then
\begin{align*}
\cosh b_2&\leq\frac 1{\sqrt{\sin^2\phi-\frac{1}{\sinh ^2 b_1}}}\leq
\Big((1-{c'_2}^2\,e^{-t})^2-
\frac{1}{\sinh ^2 (\frac{t}{2}-R-1-c_2)}\Big)^{-\frac{1}{2}}
=1+\operatorname{O}(e^{-t})\;.
\end{align*}
Since $\cosh  u\sim 1+\frac{u^2}{2}$ as $u\ra 0$, Step 3 follows.

\medskip\noindent{\bf Step 4: Conclusion. } 
Let $t\geq t_0=\max\{t_2,t_3,3\}$, $c_0=\max\{2e^2R,2(c_2+c_3)\}$ and,
with the previous notation, let $s= s^-+s^+\in\;]-2\eta,2\eta[\,$.  By
convexity, the triangle inequality and Equation
\eqref{eq:pincexpofortstab}, we have
\begin{align*}
d(p^-,p^+)&\leq d(x^-,x^+)\leq d(x^-,y)+ d(y,x^+)\\ &\leq 
(\frac{t}{2}+s^-+R\,e^{-\frac{t}{2}-s^-})+
(\frac{t}{2}+s^++R\,e^{-\frac{t}{2}-s^+})\leq t+s+c_0\,e^{-\frac{t}{2}}\;.
\end{align*}
Similarly, using Step 3 and Step 2, we have
\begin{align*}
  d(p^-,p^+)&\geq d(x^-,x^+)-d(p^-,x^-)-d(x^+,p^+)\geq
  d(x^-,x^+)-2c_3\,e^{-\frac{t}{2}}\\ &\geq  d(x^-,y)+
  d(y,x^+)-2d(y,y')-2\,c_3\,e^{-\frac{t}{2}}\\ &\geq 
(\frac{t}{2}+s^-)+
  (\frac{t}{2}+s^+) -2(c_2+c_3)\,e^{-\frac{t}{2}}\geq
  t+s-c_0\,e^{-\frac{t}{2}}\;.
\end{align*}
Let $y''$ be the closest point to $y'$ on the common perpendicular
$[p^-,p^+]$ (see the picture before this proof). Then, by Step 2, and
by convexity and Step 3, we have
$$
d(y,y'')\leq d(y,y')+d(y',y'')\leq
c_2\,e^{-\frac{t}{2}}+c_3\,e^{-\frac{t}{2}}
\leq c_0\,e^{-\frac{t}{2}}\;.
$$
This concludes the proof of Lemma \ref{lem:creationperp}.  \cqfd

\section{Counting and equidistribution of common 
perpendiculars}
\label{sec:perpendiculars}

Let $\wt M$, $x_0$, $\Ga$, $\delta$ and $M$ be as in the beginning of
Section \ref{sec:geomdyn}.

\subsection{A reminder on Patterson-Sullivan and 
skinning measures}

A family $(\mu_x)_{x\in\wt M}$ of finite measures on
$\partial_\infty\wt M$, whose support is the limit set $\Lambda\Ga$ of
$\Ga$, is a {\em Patterson-Sullivan density} for $\Ga$ if
$$
\ga_*\mu_x=\mu_{\ga x}
$$
for all $\ga\in\Ga$ and $x\in \wt M$, and if the following
Radon-Nikodym derivatives exist for all $x,y\in\wt M$ and satisfy for
(almost) all $\xi\in\partial_\infty\wt M$
$$
\frac{d\mu_x}{d\mu_y}(\xi)=e^{-\delta\,\beta_\xi(x,\,y)}\,.
$$

We fix such a family $(\mu_x)_{x\in\wt M}$.  The {\em Bowen-Margulis
  measure} on $T^1\wt M$ (associated with this Patterson-Sullivan
density) is the measure $\wtmBM$ on $T^1\wt M$ given by the density
\begin{equation}\label{eq:defigibbs}
d\wtmBM(v)=
e^{-\delta(\beta_{v_-}(\pi(v),\,x_0)\,+\,\beta_{v_+}(\pi(v),\,x_0))}\;
d\mu_{x_0}(v_-)\,d\mu_{x_0}(v_+)\,dt
\end{equation} 
in Hopf's parametrisation. The Bowen-Margulis measure $\wtmBM$ is
independent of $x_0$, and it is invariant under the actions of the
group $\Ga$ and of the geodesic flow. Thus, it defines a measure
$\mBM$ on $T^1M$ which is invariant under the quotient geodesic flow,
called the {\em Bowen-Margulis measure} on $T^1M$.  If $\mBM$ is
finite, then the Patterson-Sullivan densities are unique up to a
multiplicative constant; hence the Bowen-Margulis measure is uniquely
defined, up to a multiplicative constant. When finite and normalised
to be a probability measure, it is the unique measure of maximal
entropy of the geodesic flow, if the sectional curvature of $\wt M$
has bounded derivatives.

Babillot \cite[Thm.~1]{Babillot02b} showed that if the Bowen-Margulis
measure is finite, then it is mixing for the geodesic flow of $M$ if
the length spectrum of $M$ is not contained in a discrete subgroup of
$\RR$. This condition is satisfied, for example, if $\Ga$ has a
parabolic element, if $\Lambda\Ga$ is not totally disconnected (hence
if $M$ is compact), or if $\wt M$ is a surface or a (rank-one)
symmetric space, see for instance \cite{Dalbo99,Dalbo00}.
 
We refer to \cite{DalOtaPei00} for finiteness criteria of $\mBM$. In
particular, if $\Ga$ is geometrically finite (see for instance
\cite{Bowditch95} for a definition), if the critical exponent of the
stabiliser $\Ga_p$ of any parabolic fixed point $p$ in $\Ga$ is
strictly smaller than $\delta$, then $\mBM$ is finite.

\medskip
Let $D$ be a nonempty proper closed convex subset of $\wt M$.  The
{\em (outer) skinning measure} on $\normalout D$ (associated with the
Patterson-Sullivan density $(\mu_x)_{x\in\wt M}$) is the
measure $\wt\sigma^+_D$ on $\normalout D$
defined, using the positive endpoint homeomorphism $v\mapsto v_+$ from
$\normalout D$ to $\partial_\infty\wt M-\partial_\infty D$, by 
$$
d\wt\sigma^+_D(v)  = 
e^{-\delta\, \beta_{v_+}(P_D(v_+),\,x_0)}\;d\mu_{x_0}(v_{+}) \;,
$$
and the {\em (inner) skinning measure} on $\normalin D=\iota\normalout
D$ is the measure $\wt\sigma^-_D=\iota_*\wt\sigma^+_{D}$.  Since
$P_D(v_\pm)=\pi(v)$ for every $v\in\normalpm D$, we will often replace
$P_D(v_\pm)$ by $\pi(v)$ in the above formulas when there is no doubt
on what $v$ is.  We refer to \cite{ParPau14ETDS} for more background
and for the basic properties of these measures.

The skinning measures associated with horoballs are of particular
importance in this paper.  Let $w\in T^1\wt M$.  We denote the
skinning measures on the strong stable and strong unstable leaves
$W^{+}(w)$ and $W^{-}(w)$ of $w$ by
$$
\muss{w}=\wt \sigma^-_{H\!B_{+}(w)} \quad\quad\textrm{ and }\quad\quad
\musu{w}=\wt \sigma^+_{H\!B_-(w)}\;.
$$

\subsection{Equidistribution of endvectors of common 
perpendiculars in  $T^1\wt M$}
\label{subsec:upstairs}

Let $I$ be an index set endowed with a left action $(\ga,i)\mapsto \ga
i$ of $\Ga$.  A family $\D=(D_i)_{i\in I}$ of subsets of $\wt M$ or
$T^1\wt M$ indexed by $I$ is {\it $\Ga$-equivariant} if $\ga
D_i=D_{\ga i}$ for all $\ga\in\Ga$ and $i\in I$.  We equip the index
set $I$ with the $\Ga$-equivariant equivalence relation $\sim$ defined
by $i\sim j$ if and only if there exists $\ga\in\stab_\Ga D_i$ such
that $j=\ga i$ (or equivalently if $D_j=D_i$ and $j=\ga i$ for some
$\ga\in\Ga$).  Note that $\Ga$ acts on the left on the set of
equivalence classes $I/_\sim$.

An example of such a family is given by fixing a subset $D$ of $\wt M$
or $T^1\wt M$, by setting $I=\Ga$ with the left action by translations
$(\ga,i)\mapsto \ga i$, and by setting $D_i=i D$ for every
$i\in\Ga$. In this case, we have $i\sim j$ if and only if $i^{-1}j$
belongs to the stabiliser $\Ga_D$ of $D$ in $\Ga$, and $I/_\sim\;
=\Ga/\Ga_D$. More general examples include $\Ga$-orbits of (usually
finite) collections of subsets of $\wt M$ or $T^1\wt M$ with (usually
finite) multiplicities.

A $\Ga$-equivariant family $(A_i)_{i\in I}$ of closed subsets of $\wt
M$ or $T^1\wt M$ is  {\it locally finite} if for every
compact subset $K$ in $\wt M$ or $T^1\wt M$, the quotient set $\{i\in
I: A_i\cap K\ne\emptyset\}/_\sim$ is finite. In particular, the union
of the images of the sets $A_i$ by the map $\wt M\ra M$ or $T^1\wt
M\ra T^1M$ is closed. When $\Ga\bs I$ is finite, $(A_i)_{i\in I}$ is
locally finite if and only if, for all $i\in I$, the canonical map
from $\Ga_{A_i}\bs A_i$ to $M$ or $T^1M$ is proper, where $\Ga_{A_i}$
is the stabiliser of $A_i$ in $\Ga$.

Let $\D^-=(D^-_i)_{i\in I^-}$ and $\D^+=(D^+_j)_{j\in I^+}$ be locally
finite $\Ga$-equivariant families of nonempty proper closed convex
subsets of $\wt M$.  For every $(i,j)$ in $I^-\!\times I^+$ such that $D^-_i$
and $D^+_j$ have a common perpendicular, we denote by $\alpha_{i,\,j}$
this common perpendicular, by $\ell(\alpha_{i,\,j})$ its length, by
$v^-_{i,\,j} \in \partial^1_+ D_i^-$ its initial tangent vector and by
$v^+_{i,\,j}\in \partial^1_- D_i^+$ its terminal tangent vector. Note
that if $i'\sim i$, $j'\sim j$ and $\ga\in\Ga$, then
\begin{equation}\label{eq:equivalphvpm}
\ga\,\alpha_{i',\,j'}=\alpha_{\ga i,\,\ga j},\;\;\;\ell(\alpha_{i',\,j'})=
\ell(\alpha_{\ga i,\,\ga j})\;\;\;{\rm and}\;\;\;
\ga\, v^\pm_{i',\,j'}=v^\pm_{\ga i,\,\ga j}\;.
\end{equation}
The {\em inner and outer skinning measures} of the families $\D^\pm$ on
$T^1\wt M$ are
$$
\wt\sigma^\pm_{\D^\mp}=\sum_{i\in  I^\mp/_\sim}\wt\sigma^\pm_{D_i}\;.
$$

We will now prove that the ordered pairs of initial and terminal
tangent vectors of common perpendiculars of two locally finite
equivariant families of convex sets in $\wt M$ equidistribute towards
the product of the skinning measures of the families.

\btheo\label{theo:mainequidup} Let $\wt M$ be a complete simply
connected Riemannian manifold with pinched sectional curvature at most
$-1$.  Let $\Ga$ be a nonelementary discrete group of isometries of
$\wt M$.  Let $\D^-=(D^-_i)_{i\in I^-}$ and $\D^+=(D^+_j)_{j\in I^+}$
be locally finite $\Ga$-equivariant families of nonempty proper closed
locally convex subsets of $\wt M$. Assume that the Bowen-Margulis
measure $\mBM$ is finite and mixing for the geodesic flow. Then
$$
\lim_{t\ra+\infty} \;\delta\;\|\mBM\|\;e^{-\delta\, t}
\sum_{\substack{i\in I^-/_\sim,\; j\in I^+/_\sim, \;\ga\in\Ga\\ 
\overline{D^-_i}\cap \,\ga\overline{D^+_j}=\emptyset,\; 
\ell(\alpha_{i,\,\ga j})\leq t}} \; 
\Delta_{v^-_{i,\,\ga j}} \otimes\Delta_{v^+_{\ga^{-1}i,\,j}}\;=\; 
\wt\sigma^+_{\D^-}\otimes \wt\sigma^-_{\D^+}\;
$$
for the weak-star convergence of measures on the locally compact space
$T^1\wt M\times T^1\wt M$.  
\etheo

If $\D^-=(\ga x)_{\ga\in\Ga}$ and $\D^+=(\ga y)_{\ga\in\Ga}$ for some
$x,y\in\wt M$, this statement is a consequence of the proof of
\cite[Theo.~4.1.1]{Roblin03}.  We use the same technical initial trick
as Roblin but immediately after that we use a functional approach,
better suited to obtain error terms in Section
\ref{sect:erroterms}. We will give a reformulation in $T^1M\times
T^1M$ of this result in Section \ref{subsec:downstairs}, and some
applications to particular geometric situations in Section
\ref{sect:constcurvexamp}.

\medskip \dem We first give a scheme of the proof (see \cite[\S
8]{ParPauRev} for a more elaborate one).  The crucial observation is
that two convex subsets $D^-$ and $D^+$ have a common perpendicular of
length $t> 0$ if and only if $\flow{\frac{t}{2}} \normalout D^-$ and
$\flow{-\frac{t}{2}}\normalin D^+$ intersect.  After a reduction of
the statement, we introduce test functions $\phi^\mp_\eta$ vanishing
outside a small dynamical neighbourhood of $\normalpm D^\mp$, so that
the support of the product function
$\phi^-_\eta\circ\flow{-\frac{t}{2}} \;
\phi^+_\eta\circ\flow{\frac{t}{2}}$ detects the intersection of
$\flow{\frac{t}{2}} \normalout D^-$ and $\flow{-\frac{t}{2}}\normalin
D^+$ (using Subsection \ref{subsec:creatcomperp}). We will then use
the mixing of the geodesic flow to obtain the equidistribution result.

The estimation of the small terms occuring in the following steps 2,
4 and  5 is much more precise than what is needed to prove Theorem
\ref{theo:mainequidup}. These estimates will be useful to give a
speed of equidistribution of the initial and terminal vectors, and an
error term in the asymptotic of the counting function
$\N_{\D^-,\,\D^+}(t)$  in Theorem \ref{theo:expratecount}.

To shorten the notation, we assume from now on that the
sums as in the statement of Theorem \ref{theo:mainequidup} are for
$(i,j,\ga)$ such that $\alpha_{i,\,\ga j}$ exists, that is
$\overline{D^-_i}\cap \overline{\ga D^+_j}=\emptyset$. By Equation
\eqref{eq:equivalphvpm}, this sum is independent of the
choice of representatives of $i$ in $I^-/_\sim$ and $j$ in
$I^+/_\sim$.

\medskip\noindent{\bf Step 1: Reduction of the statement. }  By
additivity, by the local finiteness of the families $\D^\pm$, and by
the definition of $\wt \sigma^\pm_{\D^\mp}=\sum_{k\in I^\mp/_\sim}
\wt\sigma^\pm_{D^\mp_k}$, we only have to prove, for all fixed $i\in
I^-$ and $j\in I^+$, that, for the weak-star convergence of measures
on $T^1\wt M\times T^1\wt M$,
\begin{equation}\label{eq:reduconeij}
\lim_{t\ra+\infty} \;\delta\;\|\mBM\|\;e^{-\delta\, t}
\sum_{\ga\in\Ga\,:\;0<\ell(\alpha_{i,\,\ga j})\leq t}
\; \Delta_{v^-_{i,\,\ga j}}
\otimes\Delta_{v^+_{\ga^{-1}i,\,j}}\;=\; \wt\sigma^+_{D^-_i}\otimes
\wt\sigma^-_{D^+_j}\;.
\end{equation}

Let $\Omega^-$ be a Borel subset of $\normalout D^-_i$ and let
$\Omega^+$ be a Borel subset of $\normalin D^+_j$.  To simplify the
notation, let 
\begin{equation}\label{eq:notationstep1}
D^-=D^-_i, \;\;D^+=D^+_j, \;\;\alpha_\ga=\alpha_{i,\,\ga
  j}, \;\;\ell_\ga=\ell(\alpha_\ga), \;\;v^\pm_\ga=v^\pm_{i,\,\ga j}
\;\;{\rm and} \;\;\wt\sigma^\pm= \wt\sigma^\pm_{D^\mp}\,.
\end{equation}
Let $v^0_\ga$ be the tangent vector at the midpoint of $\alpha_\ga$
(see the picture below, sitting in $T^1\wt M$).

\begin{center}
\input{fig_intersection.pstex_t}
\end{center}

Assume that $\Omega^-$ and $\Omega^+$ have positive finite skinning
measures (for future use, we do not assume them to be relatively
compact), and that their boundaries in $\normalout D^-$ and $\normalin
D^+$ have zero skinning measures.  Let
$$
I_{\Omega^-,\,\Omega^+}(t)=\delta\;\|\mBM\|\;e^{-\delta \,t}\;
\card\{\ga\in\Ga\,:\;0<\ell_{\ga }\leq t, \;
v^-_{\ga}\in\Omega^-,\;v^+_{\ga}\in\ga\Omega^+\}
\;.
$$
Let us prove the stronger statement that, for every such $\Omega^\pm$,
we have
\begin{equation}\label{eq:reducnarrowup}
\lim_{t\ra+\infty} \;I_{\Omega^-,\,\Omega^+}(t)\;=\; \wt\sigma^+(\Omega^-)\;
\wt\sigma^-(\Omega^+)\;.
\end{equation}

\medskip\noindent {\bf Step 2: Construction of the bump functions. }
We recall from \cite[\S 5]{ParPau14ETDS} the definition of the test
functions $\phi^\pm_\eta$. We fix $R>0$ such that
$\mussu{w}(B^\pm(w,R)) >0$ for all $w\in \normalmp D^\pm$, hence for
all $w\in \ga\, \normalmp D^\pm$ with $\ga\in\Ga$. Such an $R$ exists
by \cite[Lem.~7]{ParPau14ETDS}.  For all $\eta,\eta'>0$, let
$h^\pm_{\eta,\,\eta'}: T^1\wt M\ra [0,+\infty[$ be the $\Ga$-invariant
measurable maps defined by
\begin{equation}\label{eq:bump}
h^\mp_{\eta,\,\eta'}(w)=
\frac{1}{2\eta\,\mussu{w}(B^\pm(w,\eta'))}
\end{equation}
if $\mussu{w}(B^\pm(w,\eta'))>0$ (which is satisfied if
$w_{\pm}\in\Lambda\Ga$), and $h^\pm_{\eta,\,\eta'}(w) =0$ otherwise.

Let us denote by $\mathbbm{1}_A$ the characteristic function
of a subset $A$. We  define the test functions $\phi^\mp_\eta=
\phi^\mp_{\eta,\,R,\,\Omega^\pm}: T^1\wt M\to[0,+\infty[$ by
\begin{equation}\label{eq:defiphi}
\phi^\mp_\eta=\;h^\mp_{\eta,\,R}\circ f^\pm_{D^\mp}\;\;
\mathbbm{1}_{\V^\pm_{\eta,\,R}(\Omega^\mp)}\,,
\end{equation}
where $\V^\pm_{\eta,\,R}(\Omega^\mp)$ and $f^\pm_{D^\mp}$ are as in
Subsection \ref{subsec:dynathick}.  Note that $v$ belongs to the
domain of definition of $f^\pm_{D^\mp}$ if $v \in \V^\pm_{\eta,\,R}
(\Omega^\mp)$, otherwise $\phi^\mp_\eta(v) = 0$. For all $v\in
T^1\wt M$ and $t\geq 0$, we have, by \cite[Lem.~17]{ParPau14ETDS},
\begin{equation}
\phi^-_{\eta,\,R,\,\Omega^+}(\flow{-t}v)= 
e\,^{-\delta\, t}\;\;
\phi^-_{\eta,\,e^{-t}R,\,\flow{t}\Omega^+}(v)\,.\label{eq:invarflottestfunct}
\end{equation}

\medskip Now, the heart of the proof is to give two pairs of upper and
lower bounds, as $T\geq 0$ is large enough and $\eta\in \;]0,1]$ is
small enough, of the quantity
\begin{equation}\label{eq:defiIetapmT}
i_{\eta}(T)=\int_0^{T}e^{\delta\,t}\;\sum_{\ga\in\Ga}\;
\int_{T^1\wt M}(\phi^-_\eta\circ\flow{-t/2})\;
(\phi^+_\eta\circ\flow{t/2}\circ\ga^{-1})\;d\wtmBM\;dt\,.
\end{equation}

\medskip
\noindent
{\bf Step 3: First upper and lower bounds. } 
For all $t\geq 0$, let
$$
a_\eta(t)=\sum_{\ga\in\Ga}\;
\int_{v\in T^1\wt M}\phi^-_\eta(\flow{-t/2}v)\;
\phi^+_\eta(\flow{t/2}\ga^{-1}v)\;d\wtmBM(v)\,.
$$
Note that by \cite[Prop.~18]{ParPau14ETDS}, we have $\int_{T^1\wt M}
\phi^\mp_\eta\;d\wtmBM=\wt \sigma^\pm(\Omega^\mp)$, which is finite
and positive. By passing to the universal cover the mixing property of
the geodesic flow on $T^1M$, for every $\epsilon>0$, there hence
exists $T_\epsilon \geq 0$ such that for all $t\geq T_\epsilon $, we
have
$$
\frac{e^{-\epsilon}}{\|\mBM\|}\;
\int_{T^1\wt M} \phi^-_\eta\;d\wtmBM\int_{T^1\wt M} \phi^+_\eta\;d\wtmBM
\leq a_\eta(t)\leq \frac{e^\epsilon}{\|\mBM\|}\;
\int_{T^1\wt M} \phi^-_\eta\;d\wtmBM\int_{T^1\wt M} \phi^+_\eta\;d\wtmBM\,.
$$
Hence for every $\epsilon>0$, there exists $c_{\epsilon,\,\eta}>0$
such that for all $T\geq 0$, we have
\begin{equation*}
e^{-\epsilon}\;
\frac{e^{\delta\,T}}{\delta\,\|\mBM\|}\;
\wt \sigma^+(\Omega^-)\;\wt \sigma^-(\Omega^+)-
c_{\epsilon,\,\eta}\leq i_{\eta}(T)\leq e^{\epsilon}\;
\frac{e^{\delta\,T}}{\delta\,\|\mBM\|}\;
\wt \sigma^+(\Omega^-)\;\wt \sigma^-(\Omega^+)+c_{\epsilon,\,\eta}\,.
\end{equation*}

\medskip
\noindent
{\bf Step 4: Second upper and lower bounds. } 
Let $T\geq 0$ and $\eta\in\;]0,1]$.  By Fubini's theorem for
nonnegative measurable maps and the definition of the test functions
$\phi^\pm_\eta$,
\begin{align}
i_{\eta}(T)=\sum_{\ga\in\Ga}\;\int_0^{T}e^{\delta\,t}\;
\int_{T^1\wt M}\;\; &h^-_{\eta,\,R}\circ f^+_{D^-}(\flow{-t/2}v)\;
h^+_{\eta,\,R}\circ f^-_{D^+}(\ga^{-1}\flow{t/2}v)\nonumber\\ &
\mathbbm{1}_{\V^+_{\eta,\,R}(\Omega^-)}(\flow{-t/2}v)\;
\mathbbm{1}_{\V^-_{\eta,\,R}(\Omega^+)}(\ga^{-1}\flow{t/2}v)
\;d\wtmBM(v)\;dt\,.\label{eq:rewriteIsetp4}
\end{align}
We start the computations by rewriting the product term involving the
technical maps $h^\pm_{\eta,\,R}$. For all $\ga\in\Ga$ and $v\in
\U^+_{D^-}\cap \U^-_{\ga D^+}$, define (using Equation
\eqref{eq:equivfibrationf})
\begin{equation}\label{eq:defwpm}
w^-=f^+_{D^-}(v)\;\;\;{\rm  and}\;\;\; 
w^+= f^-_{\ga D^+}(v)=\ga f^-_{ D^+}(\ga^{-1}v) \;.
\end{equation}
By the invariance of $f^\pm_{D^\mp}$ by precomposition by the geodesic
flow, $w^\mp$ is unchanged if $v$ is replaced by $g^sv$ for any
$s\in\RR$.  A computation, using Equation \eqref{eq:equivfibrationf}
and the $\Ga$-invariance of $h^\pm_{\eta,\,R}$, see also
\cite[p.~1334]{ParPau14ETDS}, shows that
\begin{align*}
h^-_{\eta,\,R}\circ f^+_{D^-}(\flow{-\frac{t}{2}}v)=
e^{-\delta\, \frac{t}{2}}\; 
h^-_{\eta,\,e^{-t/2}R}(\flow{t/2}w^-)\,,
\end{align*}
\begin{align*}
h^+_{\eta,\,R}\circ f^-_{D^+}(\ga^{-1}\flow{\frac{t}{2}}v)=
e^{-\delta\, \frac{t}{2}}\; 
h^+_{\eta,\,e^{-t/2}R}(\flow{-\frac{t}{2}}w^+)\;,
\end{align*}
hence, 
\begin{equation*}
h^-_{\eta,\,R}\circ f^+_{D^-}(\flow{-t/2}v)\;
h^+_{\eta,\,R}\circ f^-_{D^+}(\ga^{-1}\flow{t/2}v) =
e^{-\delta\,t}\;
h^-_{\eta,\,e^{-t/2}R}(\flow{t/2}w^-)\,
h^+_{\eta,\,e^{-t/2}R}(\flow{-t/2}w^+)\,.
\end{equation*}

The remaining product term $\mathbbm{1}_{\V^+_{\eta, \,R} (\Omega^-)}
(\flow{-t/2}v)\; \mathbbm{1}_{\V^-_{\eta,\,R}(\Omega^+)} (\ga^{-1}
\flow{t/2}v)$ in Equation \eqref{eq:rewriteIsetp4} is different from
$0$ (hence equal to $1$) if and only if
$$
v\in \flow{t/2}\V^+_{\eta, \,R}
(\Omega^-)\cap \ga\flow{-t/2} \V^-_{\eta,\,R} (\Omega^+)=
\V^+_{\eta,\,e^{-t/2}R}(\flow{t/2}\Omega^-)\cap
\V^-_{\eta,\,e^{-t/2}R}(\ga\flow{-t/2}\Omega^+)\,,
$$ 
see Section \ref{subsec:dynathick}, in particular Equation
\eqref{eq:flowbehavdynaneigh}.  By Lemma \ref{lem:creationperp}, there
exists $t_0,c_0>0$ such that for all $\eta\in \;]0,1]$ and $t\geq
t_0$, for all $v\in T^1\wt M$, if $\mathbbm{1}_{\V^+_{\eta, \,R}
  (\Omega^-)} (\flow{-t/2}v)\; \mathbbm{1}_{\V^-_{\eta,\,R}(\Omega^+)}
(\ga^{-1} \flow{t/2}v)\neq 0$, then the following facts hold:

\smallskip\noindent (i) by the convexity of $D^\pm$, we have $v\in
\U^+_{D^-}\cap \U^-_{\ga D^+}$,

\smallskip\noindent(ii) by the definition of $w^{\pm}$ (see Equation
  \eqref{eq:defwpm}), we have $w^{-}\in \Omega^{-}$ and $w^{+}\in \ga
  \Omega^{+}$ (The notation $(w^-,w^+)$ here coincides with the
  notation $(w^-, w^+)$ in Lemma \ref{lem:creationperp}),

\smallskip\noindent (iii) there exists a common perpendicular
  $\alpha_\ga$ from $D^-$ to $\ga D^+$ with $|\;\ell_\ga - t\;|\leq
  2\eta + c_0\,e^{-t/2}$, $d(\pi(v^\pm_\ga),\pi(w^\pm))\le
  c_0\,e^{-t/2}$, $d(\pi(\flow{\pm t/2}w^\mp),\pi(v))\le
  \eta+c_0\,e^{-t/2}$ and such that $\pi(v)$ is at distance at most
  $c_0\,e^{-t/2}$ from some point $p_v$ of $\alpha_\ga$.

\medskip For all $\eta\in \;]0,1]$, $\ga\in\Ga$ and $T\geq t_0$, define
$$
\A_{\eta,\ga}(T)=\big\{(t,v)\in [t_0,T] \times T^1\wt M:
v\in \V^+_{\eta,\,e^{-t/2}R}(\flow{t/2}\Omega^-)\cap
\V^-_{\eta,\,e^{-t/2}R}(\ga\flow{-t/2}\Omega^+)\big\}
$$ 
and
\begin{align}
j_{\eta,\,\ga}(T) & =
\iint_{(t,\,v)\in \A_{\eta,\ga}(T)}
  h^-_{\eta,\,e^{-t/2}R}(\flow{t/2}w^-)\;
  h^+_{\eta,\,e^{-t/2}R}(\flow{-t/2}w^+)\;dt\;d\wtmBM(v)\nonumber\\
& =\frac{1}{(2\eta)^2}
\iint_{(t,v)\in\A_{\eta,\ga}(T)} 
\frac{dt\;d\wtmBM(v)}{{\scriptstyle
\muss{w^-_t}(B^+(w^-_t,\,r_t))\;
\musu{w^+_t}(B^-(w^+_t,\,r_t))}}\,.\label{eq:adaptatJetaga}
\end{align}
with the notation
$$
r_t=e^{-t/2}R, \;\; w^-_t=\flow{t/2}w^-\;\;\;{\rm and}\;\;\;
w^+_t=\flow{-t/2}w^+\,.
$$ 

For all $s,r\in\RR$, let $\Ga_{s,r}= \{\ga\in\Ga:t_0+2+c_0\leq
\ell_{\ga} \leq s,\; v_\ga^\pm\in\N_r\Omega^\pm\}$.  By the above,
since the integral of a function is equal to the integral on any Borel
set containing its support, and since the integral of a nonnegative
function is nondecreasing in the integration domain, there hence
exists $c_4>0$ such that for all $T'\geq T\geq 0$ and $\eta\in\;]0,1]$, we
have
$$
-\,c_4+
\sum_{\ga\in\Ga_{T-\operatorname{O}(\eta+e^{-\ell_\ga/2}),
-\operatorname{O}(\eta+e^{-\ell_\ga/2})}} \!\!j_{\eta,\,\ga}(T)
\leq i_{\eta}(T) \leq
c_4+\sum_{\ga\in\Ga_{T+\operatorname{O}(\eta+e^{-\ell_\ga/2}),
\operatorname{O}(\eta+e^{-\ell_\ga/2})}} \!\!j_{\eta,\,\ga}(T')\;.
$$
We will take $T'$ to be of the form $T+ \operatorname{O}(\eta
+e^{-\ell_\ga/2})$, for a bigger $\operatorname{O}(\cdot)$ than the
one appearing in the index of the above summation.

\medskip\noindent{\bf Step 5: Conclusion. } Let $\ga\in\Ga$ be such
that $D^-$ and $\ga D^+$ have a common perpendicular with length
$\ell_\ga\geq t_0+2+c_0$. Let us prove that for all $\epsilon>0$, if
$\eta$ is small enough and $\ell_\ga$ is large enough, then for every
$T\geq \ell_\ga+ \operatorname{O}(\eta +e^{-\ell_\ga/2})$ (with the
enough's and $\operatorname{O}(\cdot)$ independent of $\ga$), we have
\begin{equation}\label{eq:step5}
 1-\epsilon\leq j_{\eta,\,\ga}(T)\leq 1+\epsilon\,.
\end{equation}

Note that $\wt\sigma^\pm(\N_\varepsilon(\Omega^\mp))$ and
$\wt\sigma^\pm (\N_{-\varepsilon}(\Omega^\mp))$ tend to $\wt\sigma^\pm
(\Omega^\mp)$ as $\varepsilon\ra 0$ (since $\wt\sigma^\pm
(\partial\Omega^\mp)=0$ as required in Step 1). Using Step 3 and Step
4, this will prove Equation \eqref{eq:reducnarrowup}, hence will
complete the proof of Theorem \ref{theo:mainequidup}.

We say that $(\wt M,\Ga)$ {\it has radius-continuous strong
  stable/unstable ball masses} if for every $\epsilon>0$, if $r\geq 1$
is close enough to $1$, then for every $v\in T^1\wt M$, if
$B^\pm(v,1)$ meets the support of $\mussu{v}$, then
$$
\mussu{v} (B^\pm(v,r))\leq e^{\epsilon} \mussu{v} (B^\pm(v,1))\,.
$$
We say that $(\wt M,\Ga)$ {\it has radius-Hölder-continuous strong
  stable/unstable ball masses} if there exists $c\in\;]0,1]$ and
$c'>0$ such that for every $\epsilon\in\;]0,1]$, if if $B^\pm(v,1)$
meets the support of $\mussu{v}$, then
$$
\mussu{v} (B^\pm(v,1+\epsilon))\leq 
e^{c'\epsilon^c} \mussu{v} (B^\pm(v,1))\,.
$$

When the sectional curvature has bounded derivatives and when $(\wt
M,\Ga)$ has radius-Hölder-continuous strong stable/unstable ball
masses, we will prove the following stronger statement: with a
constant $c_7>0$ and functions $\operatorname{O}(\cdot)$ independent
of $\ga$, for all $\eta\in\;]0,1]$ and $T\geq \ell_\ga+
\operatorname{O}(\eta + e^{-\ell_\ga/2})$, we have
\begin{equation}\label{eq:step5bis}
j_{\eta,\,\ga}(T)=
\Big(1+\operatorname{O}\Big(\frac{e^{-\ell_\ga/2}}{2\eta}\Big)\Big)^2
e^{\operatorname{O}((\eta+e^{-\ell_\ga/2})^{c_7})}\,.
\end{equation}
This stronger version will be needed for the error term estimate in
Section \ref{sect:erroterms}. In order to obtain Theorem
\ref{theo:mainequidup}, only the fact that $j_{\eta,\,\ga}(T)$ tends
to $1$ as firstly $\ell_\ga$ tends to $+\infty$, secondly $\eta$ tends
to $0$ is needed. A reader not interested in the error term
may skip many technical details below.

Let $\eta\in\;]0,1]$ and $T\geq \ell_\ga+\operatorname{O}(\eta+
e^{-\ell_\ga/2})$.  
We start the proof of Equation \eqref{eq:step5} by defining parameters
$s^+,s^-,s,v',v''$ associated with $(t,v)\in\A_{\eta,\ga}(T)$.  
\begin{center}
\input{fig_localvoga.pstex_t}
\end{center}

We have $(t,v)\in\A_{\eta,\ga}(T)$ if and only if there exist $s^\pm
\in\;]-\eta,\eta[\,$ such that
$$
g^{\mp s^\mp}v\in B^\pm(\flow{\pm t/2} w^\mp,
e^{-t/2}R)\,.
$$ 
The notation $s^\pm$ coincides with the one in the proof of Lemma
\ref{lem:creationperp} (where $(D^+,w)$ has been replaced by $(\ga
D^+, v)$). 

In order to define the parameters $s, v',v''$, we use the well known
local product structure of the unit tangent bundle in negative
curvature. If $v\in T^1M$ is close enough to $v^0_\ga$ (in particular,
$v_-\neq (v^0_\ga)_+$ and $v_+\neq (v^0_\ga)_-$), then let
$v'=f^+_{HB_-(v^0_\ga)}(v)$ be the unique element of $W^{-}(v^0_\ga)$
such that $v'_+=v_+$, let $v''=f^-_{HB_+(v^0_\ga)}(v)$ be the unique
element of $W^{+}(v^0_\ga)$ such that $v''_-=v_-$, and let $s$ be the
unique element of $\RR$ such that $g^{-s}v\in W^{+}(v')$. The map
$v\mapsto (s,v',v'')$ is a homeomorphism from a neighbourhood of
$v^0_\ga$ in $T^1\wt M$ to a neighbourhood of $(0,v^0_\ga,v^0_\ga)$ in
$\RR\times W^{-}(v^0_\ga)\times W^{+}(v^0_\ga)$. Note that if
$v=g^rv^0_\ga$ for some $r\in\RR$ close to $0$, then
$$
w^-=v^-_\ga,  \;w^+=v^+_\ga, \; s=r, \; v'=v''=v^0_\ga, \;
s^-=\frac{\ell_\ga-t}{2}+s,\;s^+=\frac{\ell_\ga-t}{2}-s\,.
$$
Up to increasing $t_0$ (which does not change Step 4, up to increasing
$c_4$), we may assume that for every $(t,v)\in \A_{\eta,\ga}(T)$, the
vector $v$ belongs to the domain of this local product structure of
$T^1\wt M$ at $v^0_\ga$. 

The vectors $v,v',v''$ are close to $v^0_\ga$ if $t$ is large and $\eta$
small, as the next result shows. We denote (also) by $d$ the
Riemannian distance induced by Sasaki's metric on $T^1\wt M$.

\blemm \label{lem:vectproche} For every $(t,v)\in \A_{\eta,\ga}(T)$, we
have $d(v,v^0_\ga), d(v',v^0_\ga), d(v'',v^0_\ga)=
\operatorname{O}(\eta+ e^{-t/2})$.
\elemm

\dem
Consider the distance $d'$ on $T^1\wt M$, defined by 
$$
\forall\;v_1,v_2\in T^1\wt M,\;\;\; d'(v_1,v_2) =
\max_{r\in[-1,0]} d\big(\pi(\flow{r}v_1),\pi(\flow{r}v_2)\big)\;.
$$
By (iii) in Step 4, we have $d(\pi(v),\pi(v^0_\ga))=
\operatorname{O}(\eta+e^{-t/2})$. By Lemma \ref{lem:minodistHam}, we
have
$$
d(\pi(\flow{-\frac{t}{2}-s^-}v),\pi(v^-_\ga))\leq 
d(\pi(\flow{-\frac{t}{2}-s^-}v),\pi(w^-))+
d(\pi(w^-),\pi(v^-_\ga))\leq 
R+c_0\,e^{-t/2}\;.
$$
By an exponential pinching argument, we hence have $d'(v,v^0_\ga)=
\operatorname{O}(\eta+ e^{-\ell_\ga/2})$.  Since $d$ and $d'$ are
equivalent (see \cite[p.~70]{Ballmann95}), we therefore
have $d(v,v^0_\ga)= \operatorname{O}(\eta+ e^{-\ell_\ga/2})$.

For all $w\in T^1\wt M$ and $V\in T_wT^1\wt M$, we may uniquely write
$V=V^{\rm su}+V^0+V^{\rm ss}$ with $V^{\rm su}\in T_wW^{-}(w)$,
$V^0\in \RR\frac{d}{dt}_{\mid t_0}\flow t w$ and $V^{\rm ss}\in
T_wW^{-}(w)$.  By \cite[\S 7.2]{PauPolSha}, Sasaki's metric (with norm
$\|\cdot\|$) is equivalent to the Riemannian metric with (product)
norm
$$
\|V\|'=\sqrt{\|\,V^{\rm su}\,\|^2+\|\,V^0\,\|^2+\|\,V^{\rm ss}\,\|^2}\,.
$$
By the dynamical local product structure of $T^1\wt M$ in the
neighbourhood of $v^0_\ga$ and by the definition of $v',v''$, the
result follows, since the exponential map of $T^1\wt M$ at $v^0_\ga$
is almost isometric close to $0$ and the projection to a factor of a
product norm is $1$-Lipschitz.  \cqfd

\medskip We use the local product structure of the Bowen-Margulis
measure to prove the following result.

\blemm \label{lem:decompgibbs}
For every $(t,v)\in \A_{\eta,\ga}(T)$, we have
$$
 dt\,d\wtmBM(v) = 
e^{\operatorname{O}(\eta+e^{-\ell_\ga/2})}\;
dt\; ds\; d\musu{v^0_\ga}(v')\;d\muss{v^0_\ga}(v'')\,.
$$
\elemm

\dem Since the above parameter $s$ differs, when $v_-,v_+$ are fixed,
only by a constant from the time parameter in Hopf's parametrisation,
we have
$$
d\wtmBM(v)  =
\frac{e^{-\delta(\beta_{v_-}(\pi(v),\,x_0)\,+\,\beta_{v_+}(\pi(v),\,x_0))}}{e^{-\delta(
\beta_{v'_{+}}(\pi(v'),\,x_0)+\beta_{v''_{-}}(\pi(v''),\,x_0))}}\;
d\musu{v^0_\ga}(v')\,d\muss{v^0_\ga}(v'')\,dt\;.
$$
As $\pi:T^1\wt M\ra \wt M$ is $1$-Lipschitz, and since $v_+=v'_+$ and
$v_-=v''_-$, the claim follows from Lemma \ref{lem:vectproche} and the
fact that the map $x\mapsto \beta_\xi(x,x_0)$ is $1$-Lipschitz for
every $\xi\in\partial_\infty\wt M$.  \cqfd

\medskip When $\ell_\ga$ is large, the submanifold $g^{\ell_\ga/2}
\Omega^-$ has a second order contact at $v^0_\ga$ with
$W^{-}(v^0_\ga)$ and similarly, $g^{-\ell_\ga/2} \Omega^+$ has a
second order contact at $v^0_\ga$ with $W^{+} (v^0_\ga)$. Let $P_\ga$
be the plane domain of $(t,s)\in\RR^2$ such that there exist
$s^\pm\in\;]-\eta,\eta[$ with $s^\mp= \frac{\ell_\ga-t}{2} \pm s +
\operatorname{O}(e^{-\ell_\ga/2})$. Note that its area is $(2\eta +
\operatorname{O}(e^{-\ell_\ga/2}))^2$.  By the above, we have (with
the obvious meaning of a double inclusion)
$$
\A_{\eta,\ga}(T)=P_\ga\times 
B^-(v^0_\ga,r_{\ell_\ga}\,e^{\operatorname{O}(\eta+e^{-\ell_\ga/2})})
\times 
B^+(v^0_\ga,r_{\ell_\ga}\,e^{\operatorname{O}(\eta+e^{-\ell_\ga/2})})\,.
$$
By Lemma \ref{lem:decompgibbs}, we hence have
\begin{multline}
\int_{\A_{\eta,\ga}(T)} dt\,d\wtmBM(v) = 
e^{\operatorname{O}(\eta+e^{-\ell_\ga/2})}\;
(2\eta + \operatorname{O}(e^{-\ell_\ga/2}))^2
\;\times  \\  
\musu{v^0_\ga}(B^-(v^0_\ga,r_{\ell_\ga}\,e^{\operatorname{O}(\eta+e^{-\ell_\ga/2})}))\;
\muss{v^0_\ga}(B^+(v^0_\ga,r_{\ell_\ga}\,e^{\operatorname{O}(\eta+e^{-\ell_\ga/2})}))\,.
\label{eq:volAetaga}
\end{multline}

The last ingredient of the proof of Step 5 is the following continuity
property of strong stable and strong unstable ball volumes as their
center varies (see \cite[Lem.~1.16]{Roblin03},
\cite[Prop.~10.16]{PauPolSha} for related properties, though we need a
more precise control for the error term in Section
\ref{sect:erroterms}).

\blemm\label{lem:variaboul} Assume that $(\wt M,\Ga)$ has
radius-continuous strong stable/unstable ball mas\-ses. There exists
$c_5>0$ such that for every $\epsilon>0$, if $\eta$ is small enough
and $\ell_\ga$ large enough, then for every $(t,v)\in \A_{\eta,\ga}(T)$,
we have
$$
\mussu{w^\mp_t}(B^\pm(w^\mp_t,r_t))=e^{\operatorname{O}(\epsilon^{c_5})}\;
\mussu{v^0_\ga}(B^\pm(v^0_\ga,r_{\ell_\ga}))\,.
$$
If we furthermore assume that the sectional curvature of $\wt M$ has
bounded derivatives and that $(\wt M,\Ga)$ has
radius-Hölder-continuous strong stable/unstable ball masses, then we
may replace $\epsilon$ by $(\eta+ e^{-\ell_\ga/2})^{c_6}$ for some
constant $c_6>0$.  
\elemm

\dem We prove the claim for $\W^{+}$, the one for $\W^{-}$ follows
similarly. The final statement is only used for the error estimates in
Section \ref{sect:erroterms}.

\medskip
\input{fig_contsssuball.pstex_t}

\medskip
Using  Equation \eqref{eq:flotVetaeta} 
and the scaling properties of skinning measure, we have
\begin{align}
\muss{w^-_t}(B^+(w^-_t,r_t))=e^{-\delta\, t/2}\muss{w^-}(B^+(w^-,R))
\label{eq:boulboul}
\end{align}
and similarly, for every $a>0$, 
\begin{equation}\label{eq:boulboulbis}
\muss{v^0_\ga}(B^+(v^0_\ga,ar_t))= e^{-\delta\, t/2}\muss{v^-_\ga}(B^+(v^-_\ga,aR))\;.
\end{equation}

Let $h^-:B^+(w^-,R)\ra W^{+}(v^-_\ga)$ be the map such that
$(h^-(v))_-=v_-$, which is well defined and a homeomorphism onto its
image if $\ell_\ga$ is large enough (since $R$ is fixed). By
\cite[Prop.~5]{ParPau14ETDS} where $C=HB_+(w^-)$,
$C'=HB_+(v^-_\ga)$, we have, for every $v\in B^+(w^-,R)$,
$$
d\muss{w^-}(v)= e^{-\delta\,\beta_{v_-}(\pi(v),\;\pi(h^-(v)))}\;d\muss{v^-_\ga}(h^-(v))\,.
$$

Let us fix $\epsilon>0$. The strong stable balls of radius $R$
centered at $w^-$ and $v^-_\ga$ are very close (see the above
picture). More precisely, recall that $R$ is fixed, and that, as seen
above, $d(\pi(w^-),\pi(v^-_\ga))= \operatorname{O}(e^{-\ell_\ga/2})$
and $d(\pi(\flow{t/2}w^-), \pi(\flow{\ell_\ga/2}v^-_\ga))=
\operatorname{O}(\eta +e^{-\ell_\ga/2})$. Therefore we have
$d(\pi(v),\pi(h^-(v)))\leq \epsilon$ for every $v\in B^+(w^-,R)$ if
$\eta$ is small enough and $\ell_\ga$ large enough. If furthermore the
sectional curvature has bounded derivatives, then by Anosov's
arguments (see for instance \cite[Theo.~7.3]{PauPolSha}) the strong
stable foliation is H\"older-continuous.  Hence we have
$d(\pi(v),\pi(h^-(v)))=\operatorname{O}((\eta+e^{-\ell_\ga/2})^{c_5})$
for every $v\in B^+(w^-,R)$, for some constant $c_5>0$, under the
additional hypothesis on the curvature.  We also have
$h^-(B^+(w^-,R))= B^+(v^-_\ga, R\,e^{\operatorname{O}(\epsilon)})$
and, under the additional hypothesis on the curvature,
$h^-(B^+(w^-,R))= B^+(v^-_\ga,
R\,e^{\operatorname{O}((\eta+e^{-\ell_\ga/2})^{c_5})})$.  Assume in
what follows that $\epsilon=(\eta+e^{-\ell_\ga/2})^{c_5}$ under the
additional hypothesis on the curvature. Since $|\beta_\xi(x,y)|\leq
d(x,y)$ for all $x,y\in\wt M$ and $\xi\in\partial_\infty\wt M$, we
hence have, for every $v\in B^+(w^-,R)$,
$$
d\muss{w^-}(v)= e^{\operatorname{O}(\epsilon)}\;
d\muss{v^-_\ga}(h^-(v))\,.
$$
The result follows by Equations \eqref{eq:boulboul},
\eqref{eq:boulboulbis} and the continuity property in the
radius. \cqfd

\medskip Now Lemma \ref{lem:variaboul} (with $\epsilon$ as in its
statement, and when its hypotheses are satisfied) implies that
\begin{align*}
&\iint_{(t,v)\in\A_{\eta,\ga}(T)} \;
\frac{dt\;d\wtmBM(v)}{\muss{w^-_t}(B^+(w^-_t,\,r_t))\;
\musu{w^+_t}(B^-(w^+_t,\,r_t))}\\
=\; & \frac{e^{\operatorname{O}(\epsilon^{c_5})}\iint_{(t,v)\in\A_{\eta,\ga}(T)} 
dt\;d\wtmBM(v)}{\muss{v^0_\ga}(B^+(v^0_\ga,\,r_t))\;
\musu{v^0_\ga}(B^-(v^0_\ga,\,r_t))}\,.
\end{align*}
By Equation \eqref{eq:adaptatJetaga} and Equation
\eqref{eq:volAetaga}, we hence have
$$
j_{\eta,\,\ga}(T)=e^{\operatorname{O}(\eta+e^{-\ell_\ga/2})}
\;e^{\operatorname{O}(\epsilon^{c_5})}\;
\frac{(2\eta + \operatorname{O}(e^{-\ell_\ga/2}))^2}{(2\eta)^2} 
$$
under the technical assumptions of Lemma \ref{lem:variaboul}. The
assumption on radius-continuity of strong stable/unstable ball
mas\-ses can be bypassed using bump functions, as explained in
\cite[p.~81]{Roblin03}, which concludes the proof of Step 5. 
\cqfd

\subsection{Equidistribution of endvectors of common 
perpendiculars in  $T^1M$}
\label{subsec:downstairs}

We now deduce from Theorem \ref{theo:mainequidup}, which is an
equidistribution result in $T^1\wt M\times T^1\wt M$, an
equidistribution result in its quotient $T^1M\times T^1M$ by the
action of $\Ga\times\Ga$.

Let $\D=(D_i)_{i\in I}$ be a locally finite $\Ga$-equivariant family
of nonempty proper closed convex subsets of $\wt M$.  Let
$\Omega=(\Omega_i)_{i\in I}$ be a $\Ga$-equivariant family of subsets
of $T^1\wt M$, where $\Omega_i$ is a measurable subset of $\normalpm
D_i$ for all $i\in I$ (the sign $\pm$ being constant).  Then
$$
\wt\sigma^\pm_{\Omega}=
\sum_{i\in I/_\sim}\wt\sigma^\pm_{D_i}|_{\Omega_i}\;,
$$ 
is a well-defined (independent of the choice of representatives in
$I/_\sim$), $\Ga$-invariant, locally finite measure on $T^1\wt M$
whose support is contained in $\bigcup_{i\in I/_\sim} \Omega_i$. The
measure $\wt\sigma^\pm_{\Omega}$ induces a locally finite measure on
$T^1M$, denoted by $\sigma^\pm_{\Omega}$.  The measures induced by
$\wt\sigma^\pm_{\D}$ on $T^1M=\Ga\bs T^1\wt M$
are called the {\em inner/outer skinning measures} of $\D$ on
$T^1M$. If $\wt x_0$ has stabiliser $\Ga_{\wt x_0}$ and
maps to $x_0\in M$, if $\D=(\ga \wt x_0)_{\ga\in\Ga}$, if $M$ has
dimension $n$, constant curvature and finite volume, then, normalising
the Patterson-Sullivan density so that $\|\mu_x\|= \Vol(\SSS^{n-1})$,
we have $\sigma^\pm_{\D}= \Vol_{T^1_{x_0}M}$ and $\|\sigma^\pm_{\D}\|=
\frac{\Vol(\SSS^{n-1})}{\card (\Ga_{\wt x_0})}$.  See Section
\ref{sect:constcurvexamp} for other examples.

Given $v\in T^1M$, we define the natural {\it multiplicity} of $v$
with respect to $\Omega$ by
$$
m_{\Omega}(v)= \frac{\card\,\{i\in I/_\sim\;:\;\wt v\in\Omega_i\}}
{\card(\stab_\Ga\wt v)}\;,
$$
for any preimage $\wt v$ of $v$ in $T^1\wt M$. The numerator and the
denominator are finite, by the local finiteness of $\D$ and the
discreteness of $\Ga$, and they depend only on the orbit $\Ga\wt v$.
The numerator takes into account the multiplicities of the images of
the elements of $\D$ in $T^1M$. Note that if $\Ga$ is torsion-free, if
$\Omega=\normalpm\D$, if for every $i\in I$ the quotient $\Ga_{D_i}\bs
D_i$ of $D_i$ by its stabiliser $\Ga_{D_i}$ maps injectively in
$M=\Ga\bs \wt M$ (by the map induced by the inclusion of $D_i$ in
$M$), and if for every $i,j\in I$ such that $j\notin \Ga i$, the
intersection $D_i\cap D_j$ is empty, then the nonzero multiplicities
$m_{\Omega}(v)$ are all equal to $1$.

Given $t> 0$ and two unit tangent vectors $v,w\in T^1M$, let
$$
n_t(v,w)=\sum_\alpha\;\;\card(\Ga_\alpha)\,,
$$
where the sum ranges over the locally geodesic paths $\alpha:[0,s]\ra
M$ such that $\dot\alpha(0)=v$, $\dot\alpha(s)=w$ and $s\in\;]0,t]$,
and $\Ga_\alpha$ is the stabiliser in $\Ga$ of any geodesic path $\wt
\alpha$ in $\wt M$ mapping to $\alpha$ by the quotient map $\wt M\ra
M$.  If $\Ga$ is torsion free, then $n_t(v,w)$ is precisely the number
of locally geodesic paths having $v$ and $w$ as initial and terminal
tangent vectors respectively, with length at most $t$.

Let $\Omega^\pm=(\Omega^\pm_i)_{i\in I^\pm}$ be $\Ga$-equivariant
families of subsets of $T^1 \wt M$, where $\Omega^\mp_k$ is a
measurable subset of $\normalpm D^\mp_k$ for all $k\in I^\mp$. We will
denote by $\N_{\Omega^-,\,\Omega^+}(t)$ the number of common
perpendiculars whose initial vectors belong to the images in $M$ of
the elements of $\Omega^-$ and terminal vectors to the images in $M$
of the elements of $\Omega^+$, counted with multiplicities:
$$
\N_{\Omega^-,\,\Omega^+}(t)=
\sum_{v,\,w\in T^1M} m_{\Omega^-}(v)\;m_{\Omega^+}(w) \;n_t(v,w)\;.
$$
When $\Omega^\pm=\normalmp \D^\pm$, we denote
$\N_{\Omega^-,\,\Omega^+}$ by $\N_{\D^-,\,\D^+}$.  If
$\Gamma$ has no torsion, if $\D^\pm=(\ga \wt D^\pm)_{\ga \in\Ga}$
where $\wt D^\pm$ is a nonempty proper closed convex subset of $\wt M$
(such that the family $\D^\pm$ is locally finite), and if $D^\pm$ is the
image of $\wt D^\pm$ by the covering map $\wt M\ra M=\Ga\bs\wt M$ (which
is a nonempty proper properly immersed closed locally convex subset of $M$),
then $\N_{\D^-,\,\D^+}$ is the counting function
$\N_{D^-,\,D^+}$ given in the introduction.

Recall that the {\it narrow topology} (also called {\em weak
  topology}) on the set $\M_{\rm f}(Y)$ of finite measures on a Polish
space $Y$ is the smallest topology such that, for every bounded
continuous map $g:Y\ra\RR$, the map from $\M_{\rm f}(Y)$ to $\RR$
defined by $\mu\mapsto \mu(g)$ is continuous.

\bcoro\label{coro:mainequicountdown} Let $\wt M,\Ga,\D^-,\D^+$
be as in Theorem \ref{theo:mainequidup}. Then,
\begin{equation}\label{eq:equidistribdown}
\lim_{t\ra+\infty}\; \delta\;\|\mBM\|\;e^{-\delta\, t}
\sum_{v,\,w\in T^1M} m_{\normalout\D^-}(v)\;m_{\normalin\D^+}(w)\; n_t(v,w) \;
\Delta_{v} \otimes\Delta_{w}\;=\;
\sigma^+_{\D^-}\otimes \sigma^-_{\D^+}\;
\end{equation}
for the weak-star convergence of measures on the locally compact space
$T^1M\times T^1M$. If $\sigma^+_{\D^-}$ and $\sigma^-_{\D^+}$ are
finite, the result also holds for the narrow convergence.

Furthermore, for all $\Ga$-equivariant families $\Omega^\pm=
(\Omega^\pm_k)_{k\in I^\pm}$ of subsets of $T^1\wt M$ with
$\Omega_k^\mp$ a Borel subset of $\partial^1_{\pm} D_k^\mp$ for all
$k\in I^\mp$, with nonzero finite skinning measure and with boundary
in $\partial^1_{\pm} D_k^\mp$ of zero skinning measure, we have as
$t\ra+\infty$
$$
\N_{\Omega^-,\,\Omega^+}(t)\;\sim\;
\frac{\|\sigma^+_{\Omega^-}\|\;\|\sigma^-_{\Omega^+}\|}{\delta\;\|\mBM\|}\;
e^{\delta\, t}\;.
$$
\ecoro

\dem Note that the sum in Equation \eqref{eq:equidistribdown}
is locally finite, hence it defines a locally finite measure on
$T^1M\times T^1M$.  We are going to rewrite the sum in the statement
of Theorem \ref{theo:mainequidup} in a way which makes it easier to
push it down from $T^1\wt M\times T^1\wt M$ to $T^1M\times T^1M$.

For every $\wt v\in T^1\wt M$, let
$$
m^\mp(\wt v)=\card\,\{k\in I^\mp/_\sim\;:\;\wt v\in\normalpm D^\mp_k\}\;,
$$
so that for every $v\in T^1M$, the multiplicity of $v$ with respect to
the family $\normalpm \D^\mp$ is 
$$
m_{\normalpm \D^\mp}(v)=\frac{m^\mp(\wt v)}{\card(\stab_\Ga\wt v)}\;,
$$
for any preimage $\wt v$ of $v$ in $T^1\wt M$.

For all $\ga\in\Ga$ and $\wt v,\wt w\in T^1\wt M$, there exists
$(i,j)\in (I^-/_\sim)\times (I^+/_\sim)$ such that $\wt v=v^-_{i,\ga
  j}$ and $\wt w=v^+_{\ga^{-1}i,j}=\ga^{-1}v^+_{i,\ga j}$ if and only
if $\ga \wt w\in\flow{\RR}\,\wt v$, there exists $i'\in I^-/_\sim$
such that $\wt v\in \normalout D^-_{i'}$ and there exists $j'\in
I^+/_\sim$ such that $\ga\wt w\in \normalin D^+_{j'}$. Then the choice
of such elements $(i,j)$, as well as $i'$ and $j'$, is free.  We hence
have
\begin{align*}
&\sum_{\substack{i\in I^-/_\sim,\; j\in I^+/_\sim, \;\ga\in\Ga\\ 
0<\ell(\alpha_{i,\,\ga j})\leq t\,,\;
v^-_{i,\,\ga j}=\wt v\,,\;v^+_{\ga^{-1}i,\,j}=\wt
w}} \; 
\Delta_{v^-_{i,\,\ga j}} \otimes\Delta_{v^+_{\ga^{-1}i,\,j}}\\\;=&
\sum_{\substack{\ga\in\Ga,\;0<s\leq t\\ 
\ga\wt w=\flow{s}\wt v}}
\card\,\big\{(i,j)\in (I^-/_\sim)\times (I^+/_\sim) \;:\;
v^-_{i,\,\ga j}=\wt v\,,\;v^+_{\ga^{-1}i,\,j}=\wt w\big\}\;
\Delta_{\wt v} \otimes\Delta_{\wt w}\\\;=&
\sum_{\substack{\ga\in\Ga,\;0<s\leq t\\ 
\ga\wt w=\flow{s}\wt v}}
 \; 
m^-(\wt v)\;m^+(\ga \wt w)\;
\Delta_{\wt v} \otimes\Delta_{\wt w}\;.
\end{align*}
Therefore
\begin{align*}
&\sum_{\substack{i\in I^-/_\sim,\; j\in I^+/_\sim, \;\ga\in\Ga\\ 
0<\ell(\alpha_{i,\,\ga j})\leq t}} \; 
\Delta_{v^-_{i,\,\ga j}} \otimes\Delta_{v^+_{\ga^{-1}i,\,j}}\\\;=&
\sum_{\wt v,\,\wt w\,\in\,T^1\wt M}\;
\card\{\ga\in\Ga:\;\exists\;s\in\;]0,t],\;
\ga\wt w=\flow{s}\wt v\}
 \;
m^-(\wt v)\;m^+(\wt w)\; \Delta_{\wt v} \otimes\Delta_{\wt w}\;.
\end{align*}

By definition, $\sigma^\pm_{\D^\mp}$ is the measure on $T^1M$ induced
by the $\Ga$-invariant measure $\wt \sigma^\pm_{\D^\mp}$, see
\cite[p.~28]{PauPolSha}. Thus the claim on weak-star convergence
follows from Theorem \ref{theo:mainequidup} and Equation
\eqref{eq:reducnarrowup} (after a similar reduction as in Step 1 of
the proof of Theorem \ref{theo:mainequidup}). Since no compactness
assumptions were made on $\Omega^\pm$ in this Step 1 in order to get
Equation \eqref{eq:reducnarrowup}, the narrow convergence claim in
Corollary \ref{coro:mainequicountdown} follows.  \cqfd
\newpage

\medskip 
Using the continuity of the pushforwards of measures for the
weak-star and the narrow topologies, applied to the basepoint maps
$\pi\times \pi$ from $T^1\wt M\times T^1\wt M$ to $\wt M\times \wt M$,
and from $T^1M\times T^1M$ to $M\times M$, we have the following
result of equidistribution of the ordered pairs of endpoints of common
perpendiculars between two equivariant families of convex sets in $\wt
M$ or two families of locally convex sets in $M$.  When $M$ has
constant curvature and finite volume, $\D^-$ is the $\Ga$-orbit of a
point and $\D^+$ is the $\Ga$-orbit of a totally geodesic cocompact
submanifold, this result is due to Herrmann \cite{Herrmann62}.

\bcoro\label{coro:mainequidbaspoint} Let $\wt
M,\Ga,\D^-,\D^+$ be as in Theorem \ref{theo:mainequidup}. Then
$$
\lim_{t\ra+\infty} \;\delta\;\|\mBM\|\;e^{-\delta\, t}\!\!
\sum_{\substack{i\in I^-/_\sim,\; j\in I^+/_\sim, \;\ga\in\Ga\\ 
0<\ell(\alpha_{i,\,\ga j})\leq t}} 
\Delta_{\pi(v^-_{i,\,\ga j})} \otimes
\Delta_{\pi(v^+_{\ga^{-1}i,\,j})}\;=\; 
\pi_*\wt\sigma^+_{\D^-}\otimes \pi_*\wt\sigma^-_{\D^+}\;,
$$
for the weak-star convergence of measures on the locally compact space
$\wt M\times \wt M$, and
\begin{align*}
&\lim_{t\ra+\infty}\; \delta\;\|\mBM\|\;e^{-\delta\, t}\!
\sum_{v,\,w\in T^1M} \!m_{\normalout\D^-}(v)\;m_{\normalin\D^+}(w)\; n_t(v,w) \;
\Delta_{\pi(v)} \otimes\Delta_{\pi(w)}\\&=\;\;
\pi_*\sigma^+_{\D^-}\otimes \pi_*\sigma^-_{\D^+}\;,
\end{align*}
for the weak-star convergence of measures on $M\times M$.  If the
measures $\sigma^\pm_{\D^\mp}$ are finite, then the above claim holds
for the narrow convergence of measures on $M\times M$. \cqfd \ecoro

Before proving Theorems \ref{theo:mainintrocount} and
\ref{theo:mainintroequidis} in the introduction, we recall the
definition of a proper nonempty properly immersed closed locally
convex subset $D^\pm$ in a negatively curved complete connected
Riemannian manifold $N$: it is a locally convex (not necessarily
connected) geodesic metric space $D^\pm$ endowed with a continuous map
$f^\pm:D^\pm\ra N$ such that, if $\wt N\ra N$ and $\wt D^\pm\ra D^\pm$
are (locally isometric) universal covers, if $\wt f^\pm:\wt
D^\pm\ra\wt N$ is a lift of $f^\pm$, then $\wt f^\pm$ is, on each
connected component of $\wt D^\pm$, an isometric embedding whose image
is a proper nonempty closed locally convex subset of $\wt N$, and the
family of images, under the covering group of $\wt N\ra N$, of the
images by $\wt f^\pm$ of the connected components of $\wt D^\pm$ is
locally finite.

\medskip\noindent{\bf Proof of Theorems \ref{theo:mainintrocount} and 
  \ref{theo:mainintroequidis}. } Let
$I^\pm=\Ga\times \pi_0(\wt D^\pm)$ with the action of $\Ga$ defined by
$\ga\cdot (\alpha,c) = (\ga \alpha,c)$ for all $\ga,\alpha\in \Ga$ and
every component $c$ of $\wt D^\pm$. Consider the families
$\D^\pm=(D^\pm_k)_{k\in I^\pm}$ where $D^\pm_k= \alpha \;\wt f^\pm(c)$
if $k=(\alpha,c)$.  Then $\D^\pm$ are $\Ga$-equivariant families of
nonempty proper closed convex subsets of $\wt M$, which are locally
finite since $D^\pm$ are properly immersed in $M$. Theorems
\ref{theo:mainintrocount} and \ref{theo:mainintroequidis} then follow
from Corollary \ref{coro:mainequicountdown}.   \cqfd

\bcoro Let $\wt M,\Ga,\D^-,\D^+$ be as in Theorem
\ref{theo:mainequidup}. Assume that $\sigma^\pm_{\D^\mp}$ are finite
and nonzero. Then
$$
\lim_{s\to+\infty}\lim_{t\to+\infty}
\frac{\delta\,\|\mBM\|^2\,e^{-\delta\, t}}
{\|\sigma^+_{\D^-}\|\|\sigma^-_{\D^+}\|} 
\sum_{v\in T^1M} m_{\normalout\D^-}(v)\;n_{t,\,\D^+}(v)\;
\Delta_{\flow s v}=\mBM\,,
$$ 
where
$$
n_{t,\D^+}(v)=\sum_{w\in T^1M}m_{\normalin\D^+}(w)\;n_t(v,w)
$$
is the number (counted with multiplicities) of locally geodesic paths
in $M$ of length at most $t$, with initial vector $v$, arriving
perpendicularly to $\D^+$. 
\ecoro

\dem For every $s\in\RR$, by Corollary \ref{coro:mainequicountdown},
using the continuity of the pushforwards of measures by the first
projection $(v,w)\mapsto v$ from $T^1M\times T^1M$ to $T^1M$, and by
the geodesic flow on $T^1M$ at time $s$, since
$(\flow{s})_*\Delta_v=\Delta_{g^sv}$, we have
$$
\lim_{t\ra+\infty}\; \delta\;\|\mBM\|\;e^{-\delta\, t}
\sum_{v\in T^1M} m_{\normalout\D^-}(v)\;n_{t,\D^+}(v) \;
\Delta_{g^sv} \;=\;
(\flow{s})_*\sigma^+_{\D^-}\|\sigma^-_{\D^+}\|\,.
$$
The result then follows from \cite[Thm.1]{ParPau14ETDS}.  \cqfd

\section{Error terms}
\label{sect:erroterms}

Let $\wt M$, $x_0$, $\Ga$, $\delta$ and $M$ be as in the beginning of
Section \ref{sec:geomdyn}.  We assume that the Bowen-Margulis measure
$\mBM$ is finite, and we define $\overline{\mBM}=
\frac{\mBM}{\|\mBM\|}$.

In this section, we give bounds for the error term in the
equidistribution and counting results of the previous section when the
geodesic flow is exponentially mixing and the (strong) stable and
unstable foliations are assumed to be at least H\"older-continuous.

There are two types of exponential mixing results available in this
context.  Firstly, when $\wt M$ is a symmetric space, then the
boundary at infinity of $\wt M$, the strong unstable, unstable,
stable, and strong stable foliations of $T^1\wt M$ are smooth. Hence
talking about leafwise $\C^\ell$ functions on $T^1M$ makes sense. Let
$\C_c^\ell(T^1M)$ be the space of $\C^\ell$ functions on $T^1M$ with
compact support and by $\|\psi\|_\ell$ the Sobolev $W^{\ell,2}$-norm
of $\psi\in\C_c^\ell(T^1M)$.

For $\ell\in \NN$, we say that the geodesic flow on $T^1M$ is
{\it exponentially mixing for the Sobolev regularity $\ell$} if there
exist $c,\kappa>0$ such that for all $\phi,\psi\in \C_c^\ell(T^1M)$
and $t\in\RR$, we have
$$
\Big|\int_{T^1M} \phi\circ g^{-t}\;\psi\;d\overline{\mBM}-\int_{T^1M}
\phi\;d\overline{\mBM}\int_{T^1M} \psi\;d\overline{\mBM}\;\Big|\leq
c\,e^{-\kappa |t|}\;\|\psi\|_\ell\;\|\phi\|_\ell\;.
$$
When $\Ga$ is an arithmetic lattice 
(the Bowen-Margulis measure then coincides, up to a multiplicative
constant, with the Liouville measure), this property, for some
$\ell\in\NN$, follows from \cite[Theorem~2.4.5]{KleMar96}, with the
help of \cite[Theorem 3.1]{Clozel03} to check its spectral gap
property, and of \cite[Lemma~3.1]{KleMar99} to deal with finite cover
problems. When $M$ has finite volume, the
conditional measures on the strong stable/unstable leaves are
homogeneous, hence $(\wt M,\Ga)$ has radius-Hölder-continuous
strong stable/unstable ball masses.

Secondly, when $\wt M$ has pinched negative sectional curvature with
bounded derivatives, then the boundary at infinity of $\wt M$, the
strong unstable, unstable, stable, and strong stable foliations of
$T^1\wt M$ are only H\"older-smooth (see for instance \cite{Brin95}
when $\wt M$ has a compact quotient and \cite[Theo.~7.3]{PauPolSha}).
Hence it is appropriate to consider H\"older functions on $T^1\wt M$.
For every $\alpha\in\;]0,1[\,$, let $\C_{\rm c} ^\alpha(X)$ be the
space of $\alpha$-H\"older-continuous real-valued functions with
compact support on a metric space $(X,d)$, endowed with the H\"older
norm
$$
\|f\|_\alpha=
\|f\|_\infty+\sup_{x,\,y\in X,\;x\neq y}\frac{|f(x)-f(y)|}{d(x,y)^\alpha}\,.
$$

For $\alpha\in\;]0,1[$, we say that the geodesic flow on $T^1M$ is
{\em exponentially mixing for the H\"older regularity $\alpha$} if
there exist $c,\kappa >0$ such that for all $\phi,\psi\in \C_{\rm
  c}^\alpha(T^1M)$ and  $t\in\RR$, we have
$$
\Big|\int_{T^1M}\phi\circ\flow{-t}\;\psi\;d\overline{\mBM}-
\int_{T^1M}\phi\; d\overline{\mBM} 
\int_{T^1M}\psi\;d\overline{\mBM}\;\Big|
\le c\;e^{-\kappa|t|}\;\|\phi\|_\alpha\;\|\psi\|_\alpha\,.
$$
This holds for compact manifolds $M$ when $M$ is locally symmetric by
\cite{Moore87}, when $M$ is two-dimensional by
\cite{Dolgopyat98}, and  when $M$ is $1/9$-pinched by
\cite[Coro.~2.7]{GiuLivPol13},  see
also \cite{MohOh15}.

\btheo\label{theo:expratecount} Let $\wt M$ be a complete simply
connected Riemannian manifold with sectional curvature at most $-1$.
Let $\Ga$ be a nonelementary discrete group of isometries of $\wt M$
and let $M=\Ga\backslash\wt M$. Assume that $(\wt M,\Ga)$ has
radius-Hölder-continuous strong stable/unstable ball masses. Let
$\D^-=(D^-_i)_{i\in I^-}$ and $\D^+=(D^+_j)_{j\in I^+}$ be locally
finite $\Ga$-equivariant families of nonempty proper closed convex
subsets of $\wt M$ with finite nonzero skinning measure.

\smallskip\noindent (1) Assume that $M$ is compact and that the
geodesic flow on $T^1M$ is mixing with exponential speed for the
H\"older regularity.  Then there exist $\alpha\in\;]0,1[$ and
$\kappa'>0$ such that for all nonnegative $\psi^\pm\in \C_{\rm
  c}^\alpha(T^1M)$, we have, as $t\ra+\infty$,
\begin{multline*}
\frac{\delta\;\|\mBM\|}{e^{\delta\, t}}
\sum_{v,\,w\in T^1M} m_{\normalout\D^-}(v)\;m_{\normalin\D^+}(w)\; n_t(v,w) \;
\psi^-(v)\,\psi^+(w)\\=
\int_{T^1M}\psi^-d\sigma^+_{\D^-} \int_{T^1M}\psi^+d\sigma^-_{\D^+}
+\operatorname{O}(e^{-\kappa't}\|\psi^-\|_\alpha\,\|\psi^+\|_\alpha)\,.
\end{multline*}

\smallskip\noindent (2) Assume that $\wt M$ is a symmetric space, that
$D^\pm_k$ has smooth boundary for every $k\in I^\pm$, that $M$ has
finite volume, and that the geodesic flow on $T^1M$ is mixing with
exponential speed for the Sobolev regularity. Then there exist
$\ell\in\NN$ and $\kappa'>0$ such that for all nonnegative maps
$\psi^\pm\in \C_{\rm c}^\ell(T^1M)$, we have, as $t\ra+\infty$,
\begin{multline*}
\frac{\delta\;\|\mBM\|}{e^{\delta\, t}}
\sum_{v,\,w\in T^1M} m_{\normalout\D^-}(v)\;m_{\normalin\D^+}(w)\; n_t(v,w) \;
\psi^-(v)\,\psi^+(w)\\=
\int_{T^1M}\psi^-d\sigma^+_{\D^-} \int_{T^1M}\psi^+d\sigma^-_{\D^+}
+\operatorname{O}(e^{-\kappa't}\|\psi^-\|_\ell\,\|\psi^+\|_\ell)\,.
\end{multline*}

Furthermore, if $\D^-$ and $\D^+$ respectively have nonzero finite
outer and inner skinning measures, if $(\wt M, \Ga)$
satisfies condition (1) or (2) above, then there exists
$\kappa''>0$ such that, as $t\ra+\infty$,
$$
\N_{\D^-,\,\D^+}(t)=
\frac{\|\sigma^+_{\D^-}\|\;\|\sigma^-_{\D^+}\|}{\delta\;\|\mBM\|}\;
e^{\delta\, t}\big(1+\operatorname{O}(e^{-\kappa'' t})\big)\;.
$$
\etheo

The maps $\operatorname{O}(\cdot)$ depend on $\wt M,\Ga,\D,$ and the
speeds of mixing.

\medskip \dem We follow the proofs of Theorem \ref{theo:mainequidup}
and Corollary \ref{coro:mainequicountdown}, adding a regularisation of
the test functions $\wt\phi^\pm_\eta$ as for the deduction of
\cite[Theo.~20]{ParPau14ETDS} from \cite[Theo.~19]{ParPau14ETDS}.

Let $\beta$ be either $\alpha\in\;]0,1]$ in the Hölder regularity case
or $\ell\in\NN$ in the Sobolev regularity case.  We fix $i\in I^-$,
$j\in I^+$, and we use the notation of Equation
\eqref{eq:notationstep1}.  Let $\wt \psi^\pm\in\C^\beta(\normalmp
D^\pm)$ be such that $\int_{T^1\wt M} \wt\psi^\pm\,
d\wt\sigma^\mp_{\D^\pm}$ is finite. Under the assumptions of Assertion
(1) or (2), we first prove the following avatar of Equation
\eqref{eq:reducnarrowup}, indicating only the required changes in its
proof: there exists $\kappa_0>0$ (independent of $\wt \psi^\pm$) such
that, as $T\ra+\infty$,
\begin{multline}
\delta\;\|\mBM\|\;e^{-\delta\, T}
\sum_{\ga\in\Ga,\,0<\ell_\ga\leq T} \; 
\wt\psi^-(v^-_\ga)\,\wt\psi^+(v^+_\ga)
\\=
\int_{\normalout D^-} \wt\psi^-\, d\wt\sigma^+
\int_{\normalin D^+} \wt\psi^+\, d\wt\sigma^-
+\operatorname{O}(e^{-\kappa_0T}\|\wt\psi^-\|_\beta\,\|\wt\psi^+\|_\beta)
\label{eq:geneth13cor19errsimp}\,.
\end{multline}

By Lemma \ref{lem:fibrationholder} and the Hölder regularity of the
strong stable and unstable foliations under the assumptions of
Assertion (1), or by the smoothness of the boundary of $D^\pm$ under
the assumptions of Assertion (2), the maps $f^\pm_{D^\mp}:
\V_{\eta,\,R}^\pm(\normalpm D^\mp)\ra \normalpm D^\mp$ are
respectively Hölder-continuous or smooth fibrations, whose fiber over
$w\in \normalpm D^\mp$ is exactly $V_{w,\,\eta,\,R}^\pm$.  By applying
leafwise the regularisation process described in the proof of
\cite[Theo.~20]{ParPau14ETDS} to characteristic functions, there
exists a constant $\kappa_1>0$ and $\chi^\pm_{\eta,\,R}\in
\C^\beta(T^1\wt M)$ such that

$\bullet$~ $\|\chi^\pm_{\eta,\,R}\|_\beta=\operatorname{O}(\eta^{-\kappa_1})$,

$\bullet$~ 
$\mathbbm{1}_{\V^\pm_{\eta \,e^{-\operatorname{O}(\eta)},
\,R\,e^{-\operatorname{O}(\eta)}}(\normalmp D^{\pm})}
\leq \chi^\pm_{\eta,\,R}\leq \mathbbm{1}_{\V^\pm_{\eta,\,R}(\normalmp D^{\pm})}$,

$\bullet$~ for every $w\in\normalmp D^{\pm}$, we have 
$$
\int_{\V^\mp_{w,\,\eta,\,R}}\chi^\pm_{\eta,\,R}\,d\nu_w^\pm=
\nu_w^\pm(\V^\mp_{w,\,\eta,\,R}) \,e^{-\operatorname{O}(\eta)}=
\nu_w^\pm(\V^\mp_{w,\,\eta \,e^{-\operatorname{O}(\eta)},\,R \,e^{-\operatorname{O}(\eta)}}) 
\,e^{\operatorname{O}(\eta)}\,,
$$
where the measures $\nu^\pm_w$ on $W^{0\mp}(v)\simeq \RR\times
W^{\mp}(v)$ are defined by $d\nu^+_w=ds\,d\musu{w}$ and
$d\nu^-_w=ds\,d\muss{w}$.  We now define the new test functions.  For
every $w\in \normalmp D^{\pm}$, let
$$
H^\pm_{\eta,\,R}(w)=
\frac{1}{\int_{\V^\mp_{w,\,\eta,\,R}}\chi^\pm_{\eta,\,R}\,d\nu_w^\pm} \,.
$$
Let $\Phi^\pm_\eta:T^1\wt M\ra \RR$ be the map defined by
$$
\Phi^\pm_\eta=(H^\pm_{\eta,\,R}\; \wt\psi^\pm)\circ f^\mp_{D^\pm}
\;\; \chi^\pm_{\eta,\,R} \,.
$$
The support of this map is contained in $\V^\pm_{\eta,\,R}(\normalmp
D^{\pm})$. Since $M$ is compact in Assertion (1) and by homogeneity in
Assertion (2), if $R$ is large enough, by the definitions of the
measures $\nu_w^\pm$, the denominator of $H^\pm_{\eta,\,R}(w)$ is a
least $c \,\eta$, where $c>0$. The map $H^\pm_{\eta,\,R}$ is hence
Hölder continuous under the assumptions of Assertion (1), and is
smooth under the ones of Assertion (2). Therefore $\Phi^\pm_\eta\in
\C^\beta(T^1\wt M)$ and there exists $\kappa_2>0$ such that
$$
\|\Phi^\pm_\eta\|_\beta=
\operatorname{O}(\eta^{-\kappa_2}\|\wt \psi^\pm\|_\beta)\,.
$$
The functions $\Phi^\mp_\eta$ are measurable, nonnegative and satisfy
$$
\int_{T^1\wt M} \Phi^\mp_\eta\;d\wtmBM=
\int_{\normalmp D^{\pm}}\wt \psi^\pm\,d\wt\sigma^\mp\,.
$$  
As in the second step of the proof of Theorem \ref{theo:mainequidup},
we will estimate in two ways the quantity
\begin{equation}\label{eq:defiIetapmTerror}
I_{\eta}(T)=\int_0^{T}e^{\delta\,t}\;\sum_{\ga\in\Ga}\;
\int_{T^1\wt M}(\Phi^-_\eta\circ\flow{-t/2})\;
(\Phi^+_\eta\circ\flow{t/2}\circ\ga^{-1})\;d\wtmBM\;dt\,.
\end{equation}

We first apply the mixing property, now with exponential decay of
correlations, as in the third step of the proof of Theorem
\ref{theo:mainequidup}. For all $t\geq 0$, let
$$
A_\eta(t)=\sum_{\ga\in\Ga}\;
\int_{v\in T^1\wt M}\Phi^-_\eta(\flow{-t/2}v)\;
\Phi^+_\eta(\flow{t/2}\ga^{-1}v)\;d\wtmBM(v)\,.
$$
Then with $\kappa>0$ as in the definitions of the exponential mixing, 
we have
\begin{align*}
A_\eta(t)&=\frac{1}{\|\mBM\|}\;
\int_{T^1\wt M}\Phi^-_\eta\,d\wtmBM\;
\int_{T^1\wt M}\Phi^+_\eta\,d\wtmBM\;+\;
\operatorname{O}(e^{-\kappa\, t}\|\Phi^-_\eta\|_\beta\|\Phi^+_\eta\|_\beta)\\ &
=\frac{1}{\|\mBM\|}\;
\int_{\normalout D^{-}}\wt \psi^-\,d\wt\sigma^+
\int_{\normalin D^{+}}\wt \psi^+\,d\wt\sigma^-\;+\;
\operatorname{O}(e^{-\kappa \,t}\eta^{-2\kappa_2}
\|\wt \psi^-\|_\beta\|\wt \psi^+\|_\beta)\,.
\end{align*}
Hence by integrating, 
\begin{equation}\label{eq:firstestimerror}
I_{\eta}(T)=
\frac{e^{\delta\,T}}{\delta\,\|\mBM\|}\;\Big(
\int_{\normalout D^{-}}\wt \psi^-\,d\wt\sigma^+
\int_{\normalin D^{+}}\wt \psi^+\,d\wt\sigma^-\;+\;
\operatorname{O}(e^{-\kappa\, T}\eta^{-2\kappa_2}
\|\wt \psi^-\|_\beta\|\wt \psi^+\|_\beta)\Big)\,.
\end{equation}

Now, we exchange the integral over $t$ and the summation over $\ga$ in
the definition of $I_{\eta}(T)$, and we proceed as in the fourth step
of the proof of Theorem \ref{theo:mainequidup}:
$$
I_{\eta}(T)=\sum_{\ga\in\Ga}\;\int_0^{T}e^{\delta\,t}\;
\int_{T^1\wt M}(\Phi^-_\eta\circ\flow{-t/2})\;
(\Phi^+_\eta\circ\flow{t/2}\circ\ga^{-1})\;d\wtmBM\;dt\,.
$$

Let $\wh\Phi^\pm_\eta=H^\pm_{\eta,\,R}\circ f^\mp_{D^\pm} \;
\chi^\pm_{\eta,\,R}$, so that $\Phi^\pm_\eta=\wt\psi^\pm\circ
f^\mp_{D^\pm}\;\wh\Phi^\pm_\eta$. By the last two properties of the
regularised maps $\chi^\pm_{\eta,\,R}$, we have, with $\phi^\mp_\eta$
defined as in Equation \eqref{eq:defiphi},
\begin{equation}\label{eq:controlwhPhi}
\phi^\pm_{\eta\,e^{-\operatorname{O}(\eta)},\,R \,e^{-\operatorname{O}(\eta),\,\normalmp D^{\pm}}}
\,e^{-\operatorname{O}(\eta)}\leq \wh\Phi^\pm_\eta\leq 
\phi^\pm_\eta\,e^{\operatorname{O}(\eta)}\;.
\end{equation}

If $v\in T^1\wt M$ belongs to the support of $(\Phi^-_\eta \circ
\flow{-t/2})\; (\Phi^+_\eta\circ \flow{t/2} \circ \ga^{-1})$, then we
have $v\in \flow{t/2}\V^+_{\eta,\,R}(\normalout D^{-})\cap \flow{-t/2}
\V^-_{\eta,\,R}(\ga\normalin D^{+})$. Hence the properties (i),
(ii) and (iii) of the fourth step of the proof of Theorem
\ref{theo:mainequidup} still hold (with $\Omega_-=\normalout D^{-}$
and $\Omega_+=\normalin (\ga D^{+})$). In particular, if
$w^-=f^+_{D^-}(v)$ and $w^+=f^-_{\ga D^+}(v)$, we have, as in the
fifth step of the proof of Theorem \ref{theo:mainequidup}, that
$$
d(w^\pm,v^\pm_\ga)=\operatorname{O}(\eta+e^{-\ell_\ga/2})\,.
$$
Hence, with $\kappa_3=\alpha$ in the Hölder case and $\kappa_3=1$ in
the Sobolev case (we may assume that $\ell\geq 1$), we have
$$
|\,\wt\psi^\pm(w^\pm)-\wt\psi^\pm(v^\pm_\ga)\,|=
\operatorname{O}((\eta+e^{-\ell_\ga/2})^{\kappa_3}\|\wt\psi^\pm\|_\beta)\,.
$$
Therefore there exists a constant $\kappa_4>0$ such that 
\begin{align*}
I_{\eta}(T)=\sum_{\ga\in\Ga}&\;(\wt\psi^-(v^-_\ga)\wt\psi^+(v^+_\ga)+
\operatorname{O}((\eta+e^{-\ell_\ga/2})^{\kappa_4}
\|\wt\psi^-\|_\beta\|\wt\psi^+\|_\beta))\times
 \\ & \int_0^{T}e^{\delta\,t}\;\int_{v\in T^1\wt M} 
\;\wh\Phi^-_{\eta}(\flow{-t/2}v)\;
\wh\Phi^+_{\eta}(\ga^{-1}\flow{t/2}v)\;d\wtmBM(v)\;dt\,.
\end{align*}

Now, using the inequalities \eqref{eq:controlwhPhi}, Equation
\eqref{eq:geneth13cor19errsimp} follows as in the last two steps of
the proof of Theorem \ref{theo:mainequidup}, by taking
$\eta=e^{-\kappa_5 T}$ for some $\kappa_5>0$.

\medskip In the same way that Corollary \ref{coro:mainequicountdown}
is deduced from Theorem \ref{theo:mainequidup}, the following result
can be deduced from Equation \eqref{eq:geneth13cor19errsimp} under the
assumptions of assertions (1) or (2).  Let $(\wt \psi_k^\pm)_{k\in
  I^\pm}$ be a $\Ga$-equivariant family of nonnegative maps
$\wt\psi_k^\pm\in\C^\beta(\normalmp D_k^\pm)$ such that $\sup_{k\in
  I^\pm} \|\wt\psi_k^\pm\|_\beta$ is finite. Extend $\wt \psi_k^\pm$
by $0$ outside $\normalmp D_k^\pm$ to define a function on $T^1\wt
M$. The measurable $\Ga$-invariant function $\wt\Psi^\pm=\sum_{k\in
  I^\pm/\sim} \wt\psi_k^\pm$ on $T^1\wt M$ defines a measurable
function $\Psi^{\,\pm}$ on $T^1M$. Assume that $\int_{T^1 M}
\Psi^{\,\pm}\, d\sigma^\mp_{\D^\pm}$ is finite. Then there exists
$\kappa'_0>0$ (independent of $(\wt \psi_k^\pm)_{k\in I^\pm}$) such
that, as $t\ra+\infty$,
\begin{align}
&\delta\;\|\mBM\|\;e^{-\delta\, t}
\sum_{v,\,w\in T^1M} m_{\normalout\D^-}(v)\;m_{\normalin\D^+}(w)\; n_t(v,w) \;
\Psi^{\,-}(v)\,\Psi^{\,+}(w)
\nonumber\\=\; &
\int_{T^1M} \Psi^{\,-}\, d\sigma^+_{\D^-}
\int_{T^1M} \Psi^{\,+}\, d\sigma^-_{\D^+}
+\operatorname{O}(e^{-\kappa'_0T}\sup_{i\in
  I^-} \|\wt\psi_i^-\|_\beta\sup_{j\in
  I^+} \|\wt\psi_j^+\|_\beta)
\label{eq:genecor19err}\,.
\end{align}
Now to prove the assertions (1) and (2), we proceed as follows. If
$\psi^\pm\in\C^\beta_c(T^1M)$, for every $k\in I^\pm$,
we denote by $\wt\psi^\pm_k$ the restriction to $\normalmp D_k^\pm$ of
$\psi^\pm\circ Tp$ where $p:\wt M\ra M$ is the universal cover. Note
that the map $\Psi^{\,\pm}$ defined above coincides with $\psi^\pm$ on
the elements $u\in T^1 M$ such that $m_{\D^\pm}(u)\neq 0$, and that
$\sup_{k\in I^\pm} \|\wt\psi_k^\pm\|_\beta \leq \|\psi^\pm\|_\beta$.
Hence the assertions (1) and (2) follow from Equation
\eqref{eq:genecor19err}.

\medskip The last statement of Theorem \ref{theo:expratecount} follows
by taking as the functions $\wt \psi_k^\pm$ the constant functions $1$
in Equation \eqref{eq:genecor19err}. \cqfd

\section{Counting closed subsets of limit sets}
\label{subsec:countinlimset}

In this section, we give counting asymptotics on very general
equivariant families of subsets of the limit sets of discrete groups of
isometries of rank one symmetric spaces.

Recall that the rank 1 symmetric spaces are the hyperbolic spaces
$\HH^n_\FF$ where $\FF$ is the set $\RR$ of real numbers, $\CC$ of
complex numbers, $\HH$ of Hamilton's quaternions, or $\OO$ of
octonions, and $n\geq 2$, with $n=2$ if $\KK=\OO$. We will normalise
them so that their maximal sectional curvature is $-1$.  We denote the
convex hull in $\HH^n_\FF$ of any subset $A$ of
$\HH^n_\FF\cup\partial_\infty\HH^n_\FF$ by $\Convhull A$.

We start with $\HH^n_\RR$. The Euclidean diameter of a subset $A$ of
the Euclidean space $\RR^{n-1}$ is denoted by $\diam A$. For any
nonempty subset $B$ of the standard sphere $\SSS^{n-1}$, we denote by
$\theta(B)$ the least upper bound of half the visual angle over pairs
of points in $B$ seen from the center of the sphere. Let $\H_\infty$
be the horoball in $\HH^n_\RR$ centred at $\infty$, consisting of the
points with vertical coordinates at least $1$. For every
Patterson-Sullivan density $(\mu_x)_{x\in\wt M}$ for a discrete
nonelementary group of isometries $\Ga$ of $\wt M$, for every horoball
$\H$ in $\wt M$, and for every geodesic ray $\rho$ starting from a
point of $\partial \H$ and converging to the point at infinity $\xi$
of $\H$, the measure $e^{\delta\, t} \mu_{\rho(t)}$ converges as $t$
tends to $+\infty$ to a measure $\mu_\H$ on $\partial_\infty \wt
M-\{\xi\}$, independent on the choice of $\rho$ (see \cite[\S
2]{HerPau04}).

\bcoro\label{coro:geneOhShaKim} Let $\Ga$ be a discrete nonelementary
group of isometries of $\HH^n_\RR$, with finite Bowen-Margulis measure
$m_{\rm BM}$.  Let $(F_i)_{i\in I}$ be a $\Ga$-equivariant family of
nonempty closed subsets in the limit set $\Lambda\Ga$, whose family
$\D^+=(\Convhull F_i)_{i\in I}$ of convex hulls in $\HH^n_\RR$ is
locally finite, with finite nonzero skinning measure.

\smallskip\noindent (1) In the upper halfspace model of $\HH^n_\RR$,
assume that $\Lambda\Ga$ is bounded in $\RR^{n-1}=
\partial_\infty \HH^n_\RR -\{\infty\}$, and that $\infty$ is not the
fixed point of an elliptic element of $\Ga$. Let $\D^-$ be the
$\Ga$-equivariant family $(\ga\H_\infty)_{\ga\in\Ga}$. Then, as
$T\ra+\infty$,
$$
\card\{i\in I/_\sim\;:\; \diam(F_i)\geq 1/T\}\sim
\;\frac{\|\sigma_{\D^-}\|\,\|\sigma_{\D^+}\|}{\delta\,\|m_{\rm BM}\|}
\;(2T)^{\delta}\;.
$$

\smallskip\noindent (2) In the unit ball model of $\HH^n_\RR$, assume
that no nontrivial element of $\Ga$ fixes $0$. As $T\to+\infty$, we
have
$$
\card\{i\in I/_\sim\;:\; \cot\theta(F_i)<T\}\sim
\frac{\|\mu_0\|\,\|\sigma_{\D^+}\|}{\delta\,\|m_{\rm
    BM}\|}\;(2T)^{\delta}\;.
$$

\smallskip\noindent (3) In the upper halfspace model of $\HH^n_\RR$,
assume that $\infty$ is not the fixed point of an elliptic element of
$\Ga$. Let $\Omega$ be a Borel subset of $\RR^{n-1}=
\partial_\infty \HH^n_\RR -\{\infty\}$ such that
$\mu_{\H_\infty}(\Omega)$ is finite and positive and
$\mu_{\H_\infty}(\partial \Omega)=0$. Then, as $T\ra+\infty$,
$$
\card\{i\in I/_\sim\;:\; \diam(F_i)\geq 1/T,\;F_i\cap
\Omega\neq\emptyset\}\sim
\;\frac{\mu_{\H_\infty}(\Omega)\,\|\sigma_{\D^+}\|}{\delta\,\|m_{\rm
    BM}\|} \;(2T)^{\delta}\;.
$$
\ecoro

This corollary generalises \cite[ Thm.~1.4]{OhSha12} and
\cite[Thm.~1.2]{OhShaCircles}) where the subsets $F_i$ are round
spheres.  When $\Ga$ is an arithmetic lattice, the error term in the
claims (1) and (2) is $\operatorname{O}(T^{\delta-\kappa})$ for some
$\kappa>0$, as it follows from Theorem \ref{theo:expratecount} (2)
 (using the Riemannian convolution smoothing process of Green and
Wu as in \cite[\S 3]{ParPau12JGA} to smooth by a very small
perturbation the boundary of $\C F_i$, so that the perturbation of the
lengths of the common perpendiculars are uniformly small).

\smallskip \dem As the Bowen-Margulis measure is finite in a locally
symmetric space, the geodesic flow is mixing (see for instance
\cite[p.~982]{Dalbo00}).

\smallskip\noindent (1) The skinning measure $\sigma^+_{\D^-}$ is
nonzero since $\Ga$ is nonelementary, and finite since the support of
$\wt \sigma^+_{\H_\infty}$, consisting of the points $v\in
\normalout\H_\infty$ such that $v_+\in\Lambda\Ga$, is compact.

For each $i\in I$, let $x_i,y_i\in F_i$ be such that $\diam
F_i=\|x_i-y_i\|$, where $\|\cdot\|$ is the Euclidean norm in
$\RR^{n-1}$.  The (signed) length $\ell(\alpha_{e,\,i})$ of the common
perpendicular $\alpha_{e,\,i}$ from $\H_\infty$ to the geodesic line
in $\HH^n_\RR$ with endpoints $x_i$ and $y_i$ (which is also the
common perpendicular from $\H_\infty$ to $\Convhull F_i$) is
$\log\frac 2{\|x_i-y_i\|}$.  Thus, since the stabiliser in $\Ga$ of an
element of $\normalout\H_\infty$ is trivial, and since $\Ga$ acts
transitively on the index set of the family $\D^-$,
\begin{align*}
\card\{i\in I/_\sim\;:\; \diam(F_i)\geq 1/T\}&=
\card\big\{i\in I/_\sim\;:\; \ell(\alpha_{e,\,i})\leq \log (2T)\big\}\\ &=
\N_{\D^-,\,\D^+}(\log (2T))\,
\end{align*}
which implies the claim (1) by Corollary \ref{coro:mainequicountdown}.

\smallskip\noindent(2) Let $\D^-$ be the $\Ga$-equivariant family
$(\{\ga 0\})_{\ga\in\Ga}$, whose skinning measure in $T^1\wt M$ is
$\wt \sigma^+_{\D^-}=\sum_{\ga\in\Ga} \mu_{\ga 0}$, so that
$\|\sigma^+_{\D^-}\|$ is equal to the (finite and nonzero) total mass
$\|\mu_0\|$ of the Patterson-Sullivan measure at $0$, since the
stabiliser of $0$ in $\Ga$ is trivial.

For each $i\in I$, let $x_i,y_i\in F_i$ be such that $\theta(F_i)=
\theta(\{x_i,y_i\})$.  The angle of parallelism formula (see for
instance \cite[p.~147]{Beardon83}) implies that $\cot\theta(F_i)=\sinh
d(0,\Convhull F_i)$, and the rest of the proof is analogous to that of
claim (1).

\smallskip\noindent (3) Note that we do not assume in (3) that the
$\Ga$-equivariant family $\D^-=(\ga\H_\infty)_{\ga\in\Ga}$ is locally
finite, and we will only use Equation \eqref{eq:reducnarrowup} (and
not Corollary \ref{coro:mainequicountdown}) to prove the claim
(3). One can check that the proof of Equation \eqref{eq:reducnarrowup}
does not use the local finiteness property of $\D^-$. By the
definition of the skinning measure $\wt \sigma^+_{\H_\infty}$ with
base point $x_0=\rho(t)$ where $\rho$ is a geodesic ray from a point
of $\partial \H_\infty$ to $\infty$, and letting $t\ra+\infty$, we see
that the pushforward of the measure $\mu_{\H_\infty}$ by the map
$x\mapsto (0,-1)\in T^1_{(x,1)}\HH^n_\RR$ from $\RR^{n-1}$ to
$\normalout \H_\infty$ is exactly the skinning measure $\wt
\sigma^+_{\H_\infty}$.  If $\diam F_i$ is small and $F_i$ meets
$\Omega$, then $F_i$ is contained in $\N_\epsilon\Omega$ for some
small $\epsilon>0$, and $\mu_{\H_\infty}(\N_\epsilon\Omega)$ converges
to $\mu_{\H_\infty}(\Omega)$ as $\epsilon\ra 0$.  We hence apply
Equation \eqref{eq:reducnarrowup} with $\Omega_e^-$ the image of
$\Omega$ by this map $x\mapsto (0,-1)$.  \cqfd

\medskip Corollary \ref{coro:geneOhShaKimintro} in the Introduction is
a special case of the following corollary. For every parabolic fixed
point $p$ of a discrete isometry group $\Ga$ of $\HH^n_\RR$, recall
that, by Bieberbach's theorem, the stabiliser of $p$ in $\Ga$ contains
a subgroup isomorphic to $\ZZ^k$ with finite index, and $k=
\operatorname{rk}_\Ga(p)\geq 1$ is called the {\it rank} of $p$ in
$\Ga$.

\bcoro\label{coro:geneOhShaKimintrogene} Let $\Ga$ be a geometrically
finite discrete group of isometries of the upper halfspace model of
$\HH^n_\RR$, whose limit set $\Lambda\Ga$ is bounded in $\RR^{n-1}
= \partial_\infty\HH^n_\RR-\{\infty\}$.  Let $\Ga_0$ be a
geometrically finite subgroup of $\Ga$ with infinite index. Assume
that the Hausdorff dimension $\delta$ of $\Lambda\Ga$ is bigger than $
\operatorname{rk}_\Ga(p)-\operatorname{rk}_{\Ga_0}(p)$ for every
parabolic fixed point $p$ of $\Ga_0$.  Then, there exists an
explicitable $c>0$ such that, as $T\ra+\infty$,
$$
\card\{\ga\in \Ga/\Ga_0\;:\; \operatorname{diam}(\ga\Lambda\Ga_0)
\geq 1/T\}\sim \;c\,T^{\delta}\;.
$$
\ecoro 

The assumption on the ranks of parabolic groups (needed to apply
\cite[Theo.~10]{ParPau14ETDS}) is in particular satisfied if every
maximal parabolic sugbroup of $\Ga_0$ has finite index in the maximal
parabolic subgroup of $\Ga$ containing it, as well as when $n=3$ and
$\delta>1$ (or equivalently if $\Ga$ does not contain a Fuchsian group
with index at most $2$, when $\Lambda\Ga$ is not totally
disconnected, see \cite[Theo.~3 (3)]{CanTay94}).

\medskip \dem First assume that $\infty$ is not fixed by an elliptic
element of $\Ga$. Since $\Ga$ is geo\-me\-tri\-cal\-ly finite, its
Bowen-Margulis measure is finite (see for instance
\cite{DalOtaPei00}).  The critical exponent of $\Ga$ is equal to the
Hausdorff dimension $\delta$ of $\Lambda\Ga$.  Let $\Ga'_0$ be the
stabiliser of the limit set $\Lambda\Ga_0$ of $\Ga_0$, and recall that
$\Ga_0$ has finite index in $\Ga'_0$ (see for instance
\cite[Coro.~4.136]{Kapovich01}).  Let us consider $I=\Ga$, the family
$(F_i=i\Lambda\Ga_0)_{i\in I}$ (which consists of nonempty closed
subsets of $\Lambda \Ga$), and $\D^+=(\C F_i)_{i\in I}$ (which is
locally finite), so that $I/_\sim= \Ga/\Ga'_0$.  Since $\Ga_0$ is
geometrically finite, the convex set $\C \Lambda \Ga_0$ is almost
cone-like in cusps and any parabolic subgroup of $\Ga$ has regular
growth (see the definitions in \cite[Sect.~4]{ParPau14ETDS}).  Hence,
under the hypothesis on the ranks of parabolic groups, by
\cite[Theo.~10]{ParPau14ETDS}, the skinning measure $\sigma^-_{\D^+}$
is finite. It is nonzero by \cite[Prop.~4 (iv)]{ParPau14ETDS}, since
$\Lambda\Ga_0\neq \Lambda\Ga$ as $\Ga_0$ has infinite index (as seen
above). Note that
\begin{multline*}
\card\{\ga\in \Ga/\Ga_0\;:\; \operatorname{diam}(\ga\Lambda\Ga_0)
\geq 1/T\}\\=[\Ga'_0:\Ga_0]\;
\card\{\ga\in \Ga/\Ga'_0\;:\; \operatorname{diam}(\ga\Lambda\Ga_0)
\geq 1/T\}\,.
\end{multline*}
The result then follows from Corollary \ref{coro:geneOhShaKim} (1).

\medskip If $\infty$ is fixed by an elliptic element of $\Ga$, let
$\Ga'$ be a finite-index torsion-free subgroup of $\Ga$ (in particular
$\Ga'$ is geometrically finite and $\Lambda\Ga'= \Lambda\Ga$).  The
action by left translations of $\Ga'$ on $\Ga/\Ga_0$ has finitely many
(pairwise distinct) orbits, say $\alpha_1\Ga_0, \dots, \alpha_k\Ga_0$.
For $i=1,\dots, k$, the group $\Ga'_i= \alpha_i\Ga_0\alpha_i^{-1}
\cap\Ga'$ is geometrically finite with infinite index in $\Ga'$. Let
$A(T)=\{\ga\Ga_0\in \Ga/\Ga_0: \operatorname{diam} (\ga\Lambda\Ga_0)
\geq 1/T\}$ and $A_i(T)= \{\ga\Ga_0\in A(T):\Ga'\ga\Ga_0=
\Ga'\alpha_i\Ga_0\}$ for $i=1,\dots, k$, so that $\card\, A(T)=
\sum_{i=1}^k \card \,A_i(T)$.  The map $\ga\Ga_0\ra \ga'\Ga'_i$ from
$A_i(T)$ to $\{\ga'\Ga'_i\in \Ga'/\Ga'_i\;:\; \operatorname{diam}
(\ga'\Lambda\Ga'_i) \geq 1/T\}$ where $\ga'\in\Ga'$ satisfies
$\ga\Ga_0=\ga'\alpha_i\Ga_0$ is easily seen to be well-defined and a
bijection. Note that the Hausdorff dimension $\delta$ of
$\Lambda\Ga'=\Lambda\Ga$ is bigger than $ \operatorname{rk}_{\Ga'}(p)-
\operatorname{rk}_{\Ga'_i}(p)= \operatorname{rk}_\Ga(\alpha_i^{-1}p)-
\operatorname{rk}_{\Ga_0} (\alpha_i^{-1}p)$ for every parabolic fixed
point $p$ of $\Ga'_i$, since $\alpha_i^{-1}p$ is a parabolic fixed
point of $\Ga_0$. By the above torsion-free case, for $i=1,\dots, k$,
there exists $c_i>0$ such that $\card \,A_i(T)\sim c_i\,T^\delta$ as
$T\ra+\infty$. The result then follows with $c=\sum_{i=1}^kc_i$.
\cqfd

\bcoro \label{coro:compdomdiscont} Let $\Ga$ be a geometrically finite
discrete subgroup of $\PSL_{2}(\CC)$ with bounded and not totally
disconnected limit set in $\CC$, which does not contain a
quasifuchsian subgroup with index at most $2$. Then there exists $c>0$
such that the number of connected components of the domain of
discontinuity $\Omega\Ga$ of $\Ga$ with diameter at least $1/T$ is
equivalent, as $T\ra+\infty$, to $c\,T^\delta$ where $\delta$ is the
Hausdorff dimension of the limit set of $\Ga$.  \ecoro

When $\infty$ is not the fixed point of an elliptic element of $\Ga$,
we have
$$
c=\frac{2^\delta\,\|\sigma_{\D^-}\|}{\delta\,\|m_{\rm BM}\|}\;
\sum_{\Omega}\|\sigma_{\wh\Omega}\|
$$ 
with $\D^-=(\ga\H_\infty)_{\ga\in\Ga}$ and $\wh\Omega=
(\ga\C\Omega)_{\ga\in\Ga}$, where $\Omega$ ranges over a set of
representatives of the $\Ga$-orbits   of the connected components
of $\Omega\Ga$ whose stabiliser has infinite index in $\Ga$.

\begin{center}
\includegraphics[trim = 0mm 15mm 0mm 15mm, clip, width=130mm]{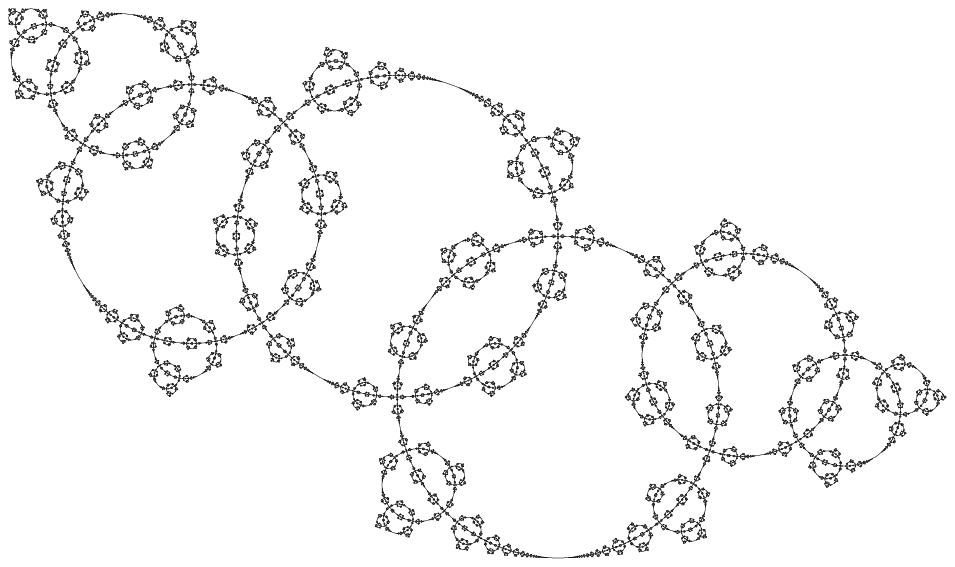}
\end{center}

\dem As mentioned after Corollary \ref{coro:geneOhShaKimintrogene}, we
have $\delta>1$, hence the assumption of this corollary on the ranks
of parabolic groups is satisfied.  By Ahlfors's finiteness theorem,
the domain of discontinuity $\Omega\Ga$ of $\Ga$ (which is a finitely
generated Kleinian group) has only finitely many orbits of connected
components (see for instance \cite[Coro.~4.108]{Kapovich01}). Since
$\Ga$ is geometrically finite, the stabiliser of a component of
$\Omega\Ga$ is again geometrically finite (see for instance
\cite[Coro.~4.112]{Kapovich01}).  

The components of $\Omega\Ga$ which are stabilised by a finite index
subgroup of $\Ga$ do not contribute to the asymptotics.  The
assumptions on $\Ga$ imply that there exists at least one other
component of $\Omega\Ga$. Otherwise indeed, the stabiliser of every
component $\Omega$ of $\Omega\Ga$ has finite index in $\Ga$, and in
particular $\partial \Omega=\Lambda\Ga$. Up to taking a finite index
subgroup, we may assume that $\Ga$ is a function group (that is,
leaves invariant a component of $\Omega \Ga$). By
\cite[Theo.~4.36]{MatTan98}, $\Ga$ is a Klein combination of
$B$-groups (groups preserving a simply connected component of their
domain of discontinuity) and elementary groups. Since $\Lambda\Ga$ is
not totally disconnected and since $\partial \Omega=\Lambda\Ga$ for
all components $\Omega$ of $\Omega\Ga$, this implies that $\Ga$ is a
$B$-group. By the structure theorem of geometrically finite $B$-groups
(see \cite[Theo.~8]{Abikoff75}), this implies that $\Ga$ is
quasifuchsian, a contradiction.  The stabilisers of these other
components of $\Omega\Ga$ have infinite index in $\Ga$. Hence the
result follows from Corollary \ref{coro:geneOhShaKimintrogene}, by a
finite summation. \cqfd

\medskip 
For example, Corollary \ref{coro:compdomdiscont} gives the
asymptotic number of the components of the domain of discontinuity
with diameter less than $\frac 1T$ as $T\to\infty$ of the crossed
Fuchsian group generated by two Fuchsian groups (see Chapter VIII \S
E.8 of \cite{Maskit88}) as in the figure below, produced using
McMullen's program { \tt lim}.  See for example Maskit's combination
theorem in op.~cit.~for a proof that crossed Fuchsian groups are
geometrically finite.

We now consider $\HH^n_\CC$, leaving to the reader the extension to
the other rank one symmetric spaces. We denote by $(w',w)\mapsto
w'\cdot w=\sum_{i=1}^{n-1} w'_i\,\overline{w_i}$ the usual Hermitian
product on $\CC^{n-1}$, and $|w|=\sqrt{w\cdot w}\,$.  Let
$$
\HH^n_\CC=\big\{(w_0,w)\in
\CC\times\CC^{n-1}\;:\; 2\,{\rm Re}\;w_0 -|w|^2>0\big\}\,,
$$
endowed with the Riemannian metric (normalised as in the beginning of
Section \ref{subsec:countinlimset})
$$
ds^2=\frac{1}{(2\,{\rm Re}\;w_0 -|w|^2)^2}
\big((dw_0-dw\cdot w)((\overline{dw_0}-w\cdot dw)+
(2\,{\rm Re}\;w_0 -|w|^2)\;dw\cdot dw\big)\,,
$$
be the Siegel domain model of the complex hyperbolic $n$-space (see
\cite[Sect.~4.1]{Goldman99}). Let 
$$
\H_\infty=\{(w_0,w)\in \CC\times\CC^{n-1}\;:\; 2\,{\rm Re}\;w_0
-|w|^2\geq 2\}\,,
$$
which is a horoball centred at $\infty$. The manifold 
$$
\Heis_{2n-1}=\partial_\infty\HH^n_\CC-\{\infty\}=\{(w_0,w)\in
\CC\times\CC^{n-1}\;:\; 2\,{\rm Re}\;w_0 -|w|^2=0\}
$$
is a Lie group (isomorphic to the $(2n-1)$-dimensional Heisenberg
group) for the law
$$
(w_0,w)\cdot (w'_0,w')=(w_0+w'_0+w\cdot w',w+w')\,.
$$
The {\it Cygan distance} $d_{\rm Cyg}$ (see \cite[p.~160]{Goldman99})
and the {\it modified Cygan distance} $d'_{\rm Cyg}$ (introduced in
\cite[Lem.~6.1]{ParPau10GT}) are the unique left-invariant distances
on $\Heis_{2n-1}$ with
$$
d_{\rm Cyg}((w_0,w),(0,0))= \sqrt{2|w_0|}\,,
\;\;\;d'_{\rm Cyg}((w_0,w),(0,0))= \sqrt{2|w_0|+|w|^2}\,.
$$
Let $d''_{\rm Cyg}=\frac{{d_{\rm Cyg}}^2}{d'_{\rm Cyg}}$, which, since
$d_{\rm Cyg}\leq d'_{\rm Cyg}\leq \sqrt{2}\,d_{\rm Cyg}$, is almost a
distance on $\Heis_{2n-1}$. For every nonempty subset $A$ of
$\Heis_{2n-1}$, we denote the ``diameter'' of $A$ for $d''_{\rm Cyg}$
by
$$
\diam_{d''_{\rm Cyg}}(A)=\max_{x,y\in
  A\,:\;x\neq y} d''_{\rm Cyg} (x,y)\;.
$$

\bcoro\label{coro:geneKim} Let $\Ga$ be a discrete nonelementary group
of isometries of the Siegel domain model of $\HH^n_\CC$, with finite
Bowen-Margulis measure $m_{\rm BM}$. Assume that $\Lambda\Ga$ is
bounded in $\Heis_{2n-1}$, and that $\infty$ is not the fixed point of
an elliptic element of $\Ga$.  Let $\D^-$ be the $\Ga$-equivariant
family $(\ga\H_\infty)_{\ga\in\Ga}$. Let $(F_i)_{i\in I}$ be a
$\Ga$-equivariant family of nonempty closed subsets in $\Lambda\Ga$,
whose family $\D^+=(\Convhull F_i)_{i\in I}$ of convex hulls in
$\HH^n_\CC$ is locally finite, with finite nonzero skinning
measure. Then, as $T\ra+\infty$,
$$
\card\{i\in I/_\sim\;:\; \diam_{d''_{\rm Cyg}}(F_i)\geq 1/T\}\sim
\frac{\|\sigma_{\D^-}\|\,\|\sigma_{\D^+}\|}{\delta\,\|m_{\rm BM}\|}
\;(2T)^{\delta}\;.
$$
\ecoro

The proof of this corollary is similar to the one of Corollary
\ref{coro:geneOhShaKim} (1), and has a similar corollary as Corollary
\ref{coro:geneOhShaKimintrogene} (replacing the rank of a parabolic
fixed point by twice the critical exponent of its stabiliser), since

$\bullet$~ the (signed) length in $\HH^n_\CC$ of the common
perpendicular from $\H_\infty$ to a geodesic in $\HH^n_\CC$ with
endpoints $x,y\in\Heis_{2n-1}$ is $\log\frac{2}{d''_{\rm Cyg}(x,y)}$
by \cite[Lem.~3.4]{ParPau11MZ};

$\bullet$~ the critical exponent of a geometrically finite group $\Ga$
of isometries of $\HH^n_\CC$ is the Hausdorff dimension of $\Lambda
\Ga$ for any (almost) distance $d_{\rm Cyg}, d'_{\rm Cyg}$ or
$d''_{\rm Cyg}$.

\section{Counting arcs in finite volume 
hyperbolic manifolds}
\label{sect:constcurvexamp}

In this subsection, we consider the special case when $M$ is a finite
volume complete connected hyperbolic good orbifold.  Taking $\wt
M=\hnr$ to be the ball model of the real hyperbolic space of dimension
$n$ and $\Ga$ to be a discrete group of isometries of $\hnr$ such that
$M$ is isometric to $\Ga\bs\hnr$, the limit set of the group $\Ga$ is
$\SSS^{n-1}$ and the Patterson-Sullivan density $(\mu_{x})_{x\in\hnr}$
of $\Ga$ can be normalised such that $\|\mu_{x}\|=\Vol(\SSS^{n-1})$
for all $x\in\hnr$.

The Bowen-Margulis measure $m_{\rm BM}$  is, by homogeneity
in this special case,  a constant multiple of the Liouville
measure $\Vol_{T^1M}$ of $T^1M$. This measure disintegrates as
$$
d\Vol_{T^1M}=\int_{x\in M} d\Vol_{T^1_xM}\;d\Vol_{M}(x)\;.
$$ 
Note that $\Vol(T^1_xM)=\frac{\Vol(\SSS^{n-1})}{\card(\Ga_{\wt x})}$
where $\Ga_{\wt x}$ is the stabiliser in $\Ga$ of any lift $\wt x$ of
$x$ in $\HH^n_\RR$, and that $\Ga_{\wt x}=\{e\}$ for $\Vol_{M}$-almost
every $x\in M$. Furthermore, if $D$ is a totally geodesic subspace or
a horoball in $\hnr$, then the skinning measures $\sigma^\pm_{D}$ are,
again by homogeneity, constant multiples of the induced Riemannian
measures $\Vol_{\normalpm D}$. These measures disintegrate with
respect to the basepoint fibration $\normalpm D\ra \partial D$ over
the Riemannian measure of the boundary $\partial D$ of $D$ in
$\HH^n_\RR$ with measure on the fiber of $x\in\partial D$ the
spherical measure on the outer/inner unit normal vectors to $D$ at
$x$:
$$
d\Vol_{\normalpm D}=
\int_{x\in \partial D} d\Vol_{\normal D\cap T^1_xM}\;
d\Vol_{\partial D}(x)\;.
$$
The following result gives the proportionality constants of the
various measures explicitly.  

\bprop\label{prop:finitevolmeasures} Let $M=\Ga\bs\hnr$ be a finite
volume orbifold of dimension $n\ge 2$. Normalise the
Patterson-Sullivan density $(\mu_{x})_{x\in\hnr}$ such that
$\|\mu_{x}\|= \Vol(\SSS^{n-1})$ for all $x\in\hnr$.
\begin{enumerate}
\item We have $m_{\rm BM}=2^{n-1}\Vol_{T^1M}$. In particular,
$$
\|m_{\rm BM}\| =2^{n-1}\Vol(\SSS^{n-1})\Vol(M)\,.
$$
\item If $D$ is a horoball in $\hnr$, then $\wt\sigma^\pm_{D}= 2^{n-1}
  \Vol_{\normalpm D}$. In particular, if $D$ is centered at a
  parabolic fixed point of $\Ga$ with stabiliser $\Ga_D$ and if
  $\D=(\ga D)_{\ga\in\Ga}$, then
$$
\|\sigma^\pm_{\D}\|=2^{n-1}\Vol(\Ga_{D}\bs\normalpm D)
=2^{n-1}(n-1)\Vol(\Ga_{D}\bs D)\;.
$$
\item If $D$ is a totally geodesic submanifold of $\hnr$ with
  dimension $k\in\{1,\dots,n-1\}$, then $\wt\sigma^+_D=\wt\sigma^-_D=
  \Vol_{\normalpm D}$. In particular, with $\Ga_D$ the stabiliser in
  $\Ga$ of $D$, if $\Ga_{D}\bs D$ is a properly immersed finite volume
  suborbifold of $M$ and if $\D=(\ga D)_{\ga\in\Ga}$, then
$$
\|\sigma^\pm_\D\|=\Vol(\Ga_{D}\bs\normalpm D)\;.
$$
If $m$ is the number of elements of $\Ga$ that pointwise
fix $D$, then
$$
\|\sigma^\pm_\D\|=\frac{1}{m}\Vol(\SSS^{n-k-1})\Vol(\Ga_{D}\bs D)\,.
$$
\end{enumerate}
\eprop

\dem Claims (1) and (3) are proven assuming that $\Ga$ has no torsion
in Proposition 10 and Claim (1) of Proposition 11 in \cite{ParPauRev},
respectively.  If $\Ga$ has torsion, Claim (1) follows by restricting
to the complement of the points in $\wt M$ with nontrivial stabiliser,
this set has zero Riemannian measure in $\wt M$, and Claim (3) follows
from the fact that the fixed point set on $D$ of an isometry which
preserves $D$, but does not pointwise fix $D$, has measure $0$ for the
Riemannian measure of $D$.

The first part of Claim (2) is proved in Claim (1) of
\cite[Prop. 10]{ParPauRev} for the outer skinning measure.  For the
second part, note that if the horoball $D$ is precisely invariant
(that is, the interiors of $D$ and $\ga D$ intersect for $\ga\in \Ga$
only if $\ga\in\Ga_{D}$), then $\Ga_{D}\bs D$ embeds in $M$ and the
image is, by definition, a Margulis cusp neighbourhood.  In the
general case, there is a precisely embedded horoball $D'$ contained in
$D$ such that $D= \N_{t}D'$ for some $t\ge 0$.  Let $\D'=(\ga
D')_{\ga\in\Ga}$. As $\Ga_{D'}=\Ga_{D}$, we have
$$
\|\sigma^\pm_{\D}\|=e^{(n-1)t}\|\sigma^\pm_{\D'}\|=
e^{(n-1)t}2^{n-1}(n-1)\Vol(\Ga_{D'}\bs D')=
2^{n-1}(n-1)\Vol(\Ga_{D}\bs D)\,,
$$
by \cite[Prop.~4 (iii)]{ParPau14ETDS}, by \cite[Prop. 10]{ParPauRev}
and  by the scaling of hyperbolic volume.  The case with torsion
follows as in Claims (1) and (3). 
\cqfd

\medskip Proposition \ref{prop:finitevolmeasures} allows us to obtain
very explicit versions of Theorems \ref{theo:mainintrocount} and
\ref{theo:4} in the case when $M$ is a finite volume hyperbolic
manifold (or good orbifold) and the properly immersed closed locally
convex subsets are any combination of points, totally geodesic
orbifolds or Margulis neighbourhoods of cusps.  The following result
gives these explicit asymptotics of the counting function in the cases
that we have not found in the literature.  We refer to the
Introduction as well as to our survey \cite{ParPauRev} for more
details and references.

\bcoro\label{coro:constcurvcount} Let $\Ga$ be a discrete group of
isometries of $\hnr$ such that $M=\Ga\bs\hnr$ has finite volume.  If
$A^-$ and $A^+$ are properly immersed finite volume totally geodesic
suborbifolds in $M$ of dimensions $k^-$ and $k^+$ in
$\{1,\dots,n-1\}$, respectively, let
$$
c(A^-,\,A^+)=\frac{\Vol(\SSS^{n-k^--1})\Vol(\SSS^{n-k^+-1})}
{2^{n-1}\,(n-1)\,\Vol(\SSS^{n-1})} 
\frac{\Vol(A^-)\Vol(A^+)}{\Vol(M)}\,.
$$
If $A^-$ and $A^+$ are Margulis cusp
neighbourhoods in $M$, let
$$
c(A^-,A^+)=\frac{2^{n-1}(n-1)\Vol(A^-)\Vol(A^+)}
{\Vol(\SSS^{n-1})\,\Vol(M)}\,.
$$
If $A^-$ is a point and $A^+$ is a Margulis cusp neighbourhood,
let
$$
c(A^-,A^+)=\frac{\Vol(A^+)}
{\Vol(M)}\,.
$$
If $A^-$ is a Margulis cusp neighbourhood and $A^+$ is a properly
immersed finite volume totally geodesic suborbifold in $M$ of
dimension $k$ in $\{1,\dots,n-1\}$, let
$$
c(A^-,A^+)=\frac{\Vol(\SSS^{n-1-k})\,\Vol(A^-)\Vol(A^+)}
{\Vol(\SSS^{n-1})\,\Vol(M)}\,.
$$
In each of these cases, if $m^\pm$ is the cardinality of the
intersection of the isotropy groups in the orbifold $M$ of the points
of $A^\pm$, then
$$
\N_{A^-,\,A^+}(t)=\N_{A^-,\,A^+,\,0}(t)\sim \frac{c(A^-,A^+)}{m^-\,m^+}\;e^{(n-1)t}\,.
$$
Furthermore, if $\Ga$ is arithmetic or if $M$ is compact, then there is
some $\kappa''>0$ such that, as $t\ra+\infty$,
$$
\N_{A^-,\,A^+}(t)=\frac{c(A^-,A^+)}{m^-\,m^+}\;
e^{(n-1)t}\big(1+\operatorname{O}(e^{-\kappa'' t})\big)\;.\;\;\;\Box
$$
\ecoro

We refer to \cite{ParPau14AFST} for several new arithmetic
applications of these results.

{\small \bibliography{../biblio} }

\bigskip
{\small\noindent \begin{tabular}{l} 
Department of Mathematics and Statistics, P.O. Box 35\\ 
40014 University of Jyv\"askyl\"a, FINLAND.\\
{\it e-mail: jouni.t.parkkonen@jyu.fi}
\end{tabular}
\medskip

\noindent \begin{tabular}{l}
D\'epartement de math\'ematique, UMR 8628 CNRS, B\^at.~425\\
Universit\'e Paris-Sud,
91405 ORSAY Cedex, FRANCE\\
{\it e-mail: frederic.paulin@math.u-psud.fr}
\end{tabular}

}

\end{document}